\documentclass[10pt]{article}
\usepackage{extarrows}
\usepackage{amsmath,amssymb,amsthm,amscd,mathrsfs}    
\usepackage{latexsym,amsmath,amssymb,graphicx}   
\usepackage[fleqn]{mathtools}
\usepackage[all]{xy} 
\usepackage{multirow}
\usepackage{color}
\usepackage{hyperref}
\usepackage{epsf}
\usepackage{tikz}
\usetikzlibrary{shapes.geometric, arrows} 
\usepackage{authblk}
\usepackage{graphicx}

\makeatletter
\newcommand\opteq[1]{\mathrel{\mathpalette\opt@eq{#1}}}
\newcommand{\opt@eq}[2]{%
  \begingroup
  \sbox\z@{$#1#2$}%
  \sbox\tw@{\resizebox{!}{.5\ht\z@}{$\m@th#1($}}%
  \nonscript\hskip-\wd\tw@
  \mkern1mu
  \raisebox{-.35\ht\z@}[0pt][0pt]{\resizebox{!}{.5\ht\z@}{$\m@th#1($}}%
  \mkern-1mu
  {#2}%
  \mkern-1mu
  \raisebox{-.35\ht\z@}[0pt][0pt]{\resizebox{!}{.5\ht\z@}{$\m@th#1)$}}%
  \mkern1mu
  \nonscript\hskip-\wd\tw@
  \endgroup
}
\makeatother

\xyoption{v2} \xyoption{2cell} \xyoption{dvips}
\CompileMatrices \numberwithin{equation}{section}
\newtheorem{prop}{Proposition}[section]
\newtheorem{theo}[prop]{Theorem}
\newtheorem{lemm}[prop]{Lemma}

\newtheorem{exam}[prop]{Example}

\newtheorem{conj}[prop]{Conjecture}

\numberwithin{equation}{section}
\makeatletter
\newcommand{\subalign}[1]{%
  \vcenter{%
    \Let@ \restore@math@cr \default@tag
    \baselineskip\fontdimen10 \scriptfont\tw@
    \advance\baselineskip\fontdimen12 \scriptfont\tw@
    \lineskip\thr@@\fontdimen8 \scriptfont\thr@@
    \lineskiplimit\lineskip
    \ialign{\hfil$\m@th\scriptstyle##$&$\m@th\scriptstyle{}##$\crcr
      #1\crcr
    }%
  }
}
\makeatother
\newcommand{\be}{\begin{equation}}
\newcommand{\ee}{\end{equation}}
\newcommand{\IP}{\mathbb{P}}%{{\relax{\rm I\kern-.18em P}}}

\newcommand\IZ{\mathbb {Z}}
\newcommand{\CM}{{\mathcal M}}

\newcommand{\umu}{{\underline \mu}}

\newcommand{\IC}{\mathbb{C}}

\newcommand{\IR}{\mathbb{R}}

\newcommand{\ba}{\begin{array}}
\newcommand{\ea}{\end{array}}

\newcommand{\tr}{{\hbox{Tr}}}

\newcommand{\CW}{{\mathcal W}}

\newcommand{\bal}{\begin{aligned}}
\newcommand{\eal}{\end{aligned}}

\newcommand{\mfg}{{\mathfrak g}}

\newcommand{\mfh}{{\mathfrak h}}

\newcommand{\wq}{{\widetilde {q}}}

\newcommand{\cW}{{\mathcal W}}
\newcommand{\cA}{{\mathcal A}}
\newcommand{\CA}{{\mathcal A}}
\newcommand{\mgg}{{\mathfrak g}}
\newcommand{\gl}{{\mathfrak l}}
\newcommand{\gs}{{\mathfrak s}}
\newcommand{\ga}{{\mathfrak a}}
\newcommand{\gb}{{\mathfrak b}}
\newcommand{\KL}{{\rm{KL}}}

\newcommand{\ch}{{\mathrm{ch}}}

\newcommand{\bfmu}{{\boldsymbol \mu}}

\CompileMatrices \LaTeXdiagrams \UseAllTwocells
\newdimen\tableauside\tableauside=1.0ex
\newdimen\tableaurule\tableaurule=0.4pt
\newdimen\tableaustep
\def\phantomhrule#1{\hbox{\vbox to0pt{\hrule height\tableaurule width#1\vss}}}
\def\phantomvrule#1{\vbox{\hbox to0pt{\vrule width\tableaurule height#1\hss}}}
\def\sqr{\vbox{%
  \phantomhrule\tableaustep
  \hbox{\phantomvrule\tableaustep\kern\tableaustep\phantomvrule\tableaustep}%
  \hbox{\vbox{\phantomhrule\tableauside}\kern-\tableaurule}}}
\def\squares#1{\hbox{\count0=#1\noindent\loop\sqr
  \advance\count0 by-1 \ifnum\count0>0\repeat}}
\def\tableau#1{\vcenter{\offinterlineskip
  \tableaustep=\tableauside\advance\tableaustep by-\tableaurule
  \kern\normallineskip\hbox
    {\kern\normallineskip\vbox
      {\gettableau#1 0 }%
     \kern\normallineskip\kern\tableaurule}%
  \kern\normallineskip\kern\tableaurule}}
\def\gettableau#1 {\ifnum#1=0\let\next=\null\else
  \squares{#1}\let\next=\gettableau\fi\next}
\allowdisplaybreaks 

\DeclareMathAlphabet{\mathpzc}{OT1}{pzc}{m}{it}
\begin{document}

\title{Affine Laumon spaces and iterated $\CW$-algebras} 
\author{Thomas Creutzig, Duiliu-Emanuel Diaconescu, Mingyang Ma}
\newcommand{\addresses}{
\footnotesize
Thomas Creutzig, \textsc{Department of Mathematics and Statistics, University of Alberta, Canada}\par\nopagebreak\textit{E-mail address}:  \texttt{creutzig@ualberta.ca}

\medskip 

Duiliu-Emanuel Diaconescu, 
\textsc{New High Energy Theory Center, Rutgers University, USA}
\par\nopagebreak\textit{E-mail address}:  \texttt{duiliu@physics.rutgers.edu}

\medskip 

Mingyang Ma, \textsc{Department of Physics and Astronomy, Rutgers University, USA} 
\par\nopagebreak\textit{E-mail address}:  \texttt{mm2595@scarletmail.rutgers.edu} }

\date{}
\maketitle

\begin{abstract} 
A family of vertex algebras whose  universal 
Verma modules coincide with the cohomology of affine Laumon spaces is found. 
This result is based 
on an explicit expression for the generating function of Poincar\'e polynomials of these spaces. 

There is  a variant of quantum Hamiltonian reduction that realizes vertex algebras which we call iterated $\CW$-algebras and our main conjecture is that the vertex algebras associated to the affine Laumon spaces are subalgebras of 
iterated $\CW$-algebras.
\end{abstract}

%\begin{abstract} 
%A conjectural construction is presented of  
%a vertex operator algebra whose universal 
%Verma module coincides with the cohomology of the affine Laumon space,
%subject to certain conditions on the numerical type. 
%This conjecture is based 
%on an explicit expression for the generating function of Poincar\'e polynomials of these spaces. The associated vertex algebra is constructed as a subalgebra of an 
%iterated $\CW$-algebra.
%\end{abstract}

\tableofcontents

\section{Introduction}

Affine Laumon spaces provide a rigorous mathematical construction for the surface defects introduced by Gukov and Witten in \cite{Ramified_L}. In the context of AGT correspondence
\cite{AGT}, this predicts a vertex operator algebra action on their equivariant 
Borel-Moore homology. 
 As shown by Schiffmann and Vasserot \cite{W_instantons} and by 
 Maulik and Okounkov \cite{Q_groups},
 for the moduli space of rank $r$ unramified instantons, this vertex algebra is the 
 $\CW$-algebra associated to a principal nilpotent element in ${\mathfrak gl}_r$. 
 More precisely, the equivariant Borel-Moore homology of the instanton moduli space is identified to a highest weight $\CW({\mathfrak gl}_r)$-module. This action is constructed using an algebraic 
 relation between the $\CW$-algebra and the Yangian of affine ${\mathfrak gl}_1$. 
 The latter is then  shown to act on the equivariant Borel-Moore homology via Hecke 
 operators. The analogous construction of moduli spaces of $G$-instantons was 
 carried out by Braverman, Fikelberg and Nakajima \cite{Inst_W}. 
 Earlier results on the representation theoretic structure of the equivariant Borel-Moore homology 
 of affine Laumon spaces of type $(1, \ldots, 1)$ include \cite{Y_rings,AL_conj,AL_int}.

A similar program for ramified instantons was initiated by Finkelberg and Tsymbaliuk  \cite{M_slices} and Negut \cite{AGT_parabolic}. Using Hecke transformations, they constructed a 
shifted affine Yangian action on the equivariant Borel-Moore homology of affine Laumon 
spaces, respectively a shifted quantum toroidal action on their equivariant K-theory. 
The main missing piece of the puzzle is the construction of the associated 
vertex operator algebras. This is understood only for the rectangular case i.e. affine Laumon spaces of type $(s, \ldots, s)$, 
which is treated in \cite{Def_W}.  Similar relations between affine Yangians and rectangular 
$\CW$-algebras were established in \cite{Super_rectangular,Rectangular_D}. 

In the physics literature, the relation between the Nekrasov partition functions 
of affine Laumon spaces and $\CW$-characters was studied in 
\cite{Kozcaz:2010af,Wyllard:2010rp,Wyllard:2010vi,Kanno:2011fw}. Most of these results are restricted to the case 
to affine Laumon spaces of type $(1, \ldots, 1)$, with the exception of 
 Kanno and Tachikawa \cite{Kanno:2011fw} who address some aspects 
of the partition function for Euler characteristics in the general case.

The present paper presents a conjectural construction of 
vertex operator algebras associated to affine Laumon spaces 
of type ${\bf r}= (r_0, \ldots, r_{\ell-1})$, $\ell\geq 2$, subject to the conditions 
$r_0 \geq r_1 \geq \ldots \geq r_{\ell-1}\geq 1$. The main result identifies the 
generating function 
for  Poincar\'e polynomials of affine Laumon spaces of type ${\bf r}$ as above to the conformal weight refinement of the universal Verma module character of the associated 
vertex algebra. These vertex operator algebras are  obtained as natural subalgebras of iterated $\CW$-algebras. 

In more detail, recall that the affine Laumon space of type $(r_0, \ldots, r_\ell)$ is 
a moduli space of framed parabolic torsion free sheaves on $\IP^1\times \IP^1$. 
The parabolic structure is specified by a flag  of subsheaves of ranks
$r_0, \ldots, r_\ell$ along 
$\IP^1\times \{0\}$, while the framing condition 
 rigidifies the parabolic sheaves along $\IP^1\times \{\infty\}\, \cup\, \{\infty\}\times \IP^1$. These spaces were studied in detail in \cite{QZ}. In particular, as shown in 
 loc. cit,  each such space is naturally identified to a connected component 
 of the fixed locus for a certain ${\bfmu}_\ell$-action on the moduli space of rank 
 $r=r_0+\cdots + r_\ell$ instantons. Moreover, using the ADHM construction, 
 the affine Laumon spaces are identified to moduli spaces of cyclic framed 
 representations of the chain-saw quiver.  Section \ref{ALspaces} provides a summary 
 of all these results. 
 
 For fixed 
 type ${\bf r}=(r_0, \ldots, r_\ell)$, the moduli space of framed parabolic sheaves 
 consists of infinitely many smooth components $\CM({\bf r},{\bf n})$. The index ${\bf n}\in 
 \IZ_{\geq 0}^{\times \ell}$ encodes the parabolic second Chern class. 
 The generating function 
 \[
Z_{\bf r } = \sum_{{\bf n}} P_y(\CM({\bf r}, {\bf n})) \prod_{a=0}^{\ell-1} q_a^{n_a} 
\]
for Poincar\'e polynomials is computed in Section \ref{PPsection} using Atiyah-Bott-Morse 
localization. This computation generalizes the results obtained by Nakajima and Yoshioka 
\cite{Lect_inst}  for unramified instantons, and employs similar techniques. Similar results were obtained by  Bruzzo et al in \cite{Framed_H} for 
framed torsion free sheaves on Hirzebruch surfaces. 

In order to present the result in a concise form, for any $a\in \IZ$, let $q_a=q_c$ and 
$r_a = r_c$ where $c\in \{0, \ldots, \ell-1\}$ is uniquely defined by $a\equiv c$ mod $\ell$.
Then one has: 

\begin{theo}\label{mainthmA} ({\it Theorem \ref{PPthmA}})
The following identity holds
\be\label{eq:introA}
Z_{\bf r} =  
\prod_{a=0}^{\ell-1} \prod_{t=1}^{r_a} \prod_{c=1}^{\ell-1}
\prod_{k=1}^{\infty} {1\over 1- y^{-2t} z^{k} }
{1\over 
1- y^{-2t} z^{k}(\wq_{a(\beta)-c})^{-1} (\wq_{a(\beta)-c-1})^{-1} \cdots 
(\wq_{a(\beta)-\ell+1})^{-1}}~,
\ee
where 
\be\label{eq:introB}
\wq_a = y^{2r_a} q_a
\ee
for any $a\in \IZ$ and 
 \[
 z=\prod_{a=0}^{\ell-1} \wq_a.
\]
 \end{theo} 
\noindent 
Moreover, note that the odd degree cohomology of the spaces $\CM({\bf r}, {\bf n})$ is trivial since the fixed locus of the torus action is finite for any ${\bf r}$, ${\bf n}$.

As shown in Proposition \eqref{newvarprop}, an alternative expression for the partition function \eqref{eq:introA}  is
\be\label{eq:introC}
\bal 
  Z_{\bf r} =  & \bigg(\prod_{a=0}^{\ell-1} \prod_{t=1}^{r_a} \prod_{k=1}^\infty {1\over 1-y^{-2t}z^k}  \bigg)  \\
  &
\prod_{1\leq a< c\leq \ell} \bigg( 
\prod_{t=1}^{r_{\ell-c}} \prod_{k=1}^\infty
{1\over 1- y^{-2t} z^{k} u_{a}u_{c}^{-1}}\bigg)
\bigg( \prod_{t=1}^{r_{\ell-a}}\prod_{k=1}^\infty
{1\over 1- y^{-2t} z^{k-1} u_{a}^{-1}u_{c}}\bigg),\\
\eal 
 \ee
 where the variables $u_c$, $1\leq c\leq \ell$, are given by 
\be\label{eq:introD}
 u_c  = (\wq_{-c})^{-1} \cdots 
(\wq_{-\ell+1})^{-1}~, \qquad 1\leq c\leq \ell-1, \qquad u_\ell=1. 
\ee

In section \ref{section:walgebras} we find a family of vertex 
algebras indexed by ${\bf r} = (r_0, \ldots, r_{\ell-1})$, $r_0\geq \cdots \geq r_{\ell-1}\geq 1$, so that the partition function \eqref{eq:introC} coincides with the conformal weight 
refinement of the character of the associated universal Verma module. These vertex algebras are $\CW$-algebras times certain commutative vertex algebras.
A $\CW$-algebra is constructed from an affine vertex algebra via quantum Hamiltonian reduction, i.e. as a BRST cohomology associated to a nilpotent element. We are interested in a variant or better iteration of quantum Hamiltonian reduction that realizes this family of vertex algebras. 

 Most of section \ref{section:walgebras} is thus  devoted to the explanation of a concept that we call iterated quantum Hamiltonian reduction. This procedure is  related to the gluing of corner vertex algebras  \cite{Gaiotto:2017euk, Creutzig:2017uxh, Prochazka:2017qum}.
Our main conjecture is then that the families of vertex algebras that are associated to affine Laumon spaces are subalgebras of iterated quantum Hamiltonian reductions.

 In order to set the stage, sections \ref{section:affine} -- \ref{characters} provide a review of $\CW$-algebras focusing on $\CW$-algebras associated to nilpotent elements in ${\mathfrak gl}_N$. In particular Sections \ref{properties} and \ref{characters} summarize their main properties as well as 
the main properties of their universal Verma modules and characters. 

The  construction of iterated $\CW$-algebras and their universal Verma modules,
which is the main construction in this paper, is presented in Section \ref{iteratedW}. 
The input data consists of: 
\begin{itemize} 
\item a simple Lie algebra ${\mathfrak g}$,
\item  an ordered sequence $({\mathfrak a}_1,\ldots, {\mathfrak a}_L)$ of subalgebras of 
${\mathfrak g}$ so that ${\mathfrak a}_1\oplus \cdots \oplus {\mathfrak a}_L$ is a 
subalgebra of ${\mathfrak g}$, and 
\item a sequence $(f_1, \ldots, f_L)$ of nilpotent elements of $({\mathfrak a}_1,\ldots, {\mathfrak a}_L)$ 
\end{itemize} 
Given the above data, 
Section \ref{iteratedW} constructs a $\CW$-algebra 
$\CW^k_L({\mathfrak g},({\mathfrak a}_1,\ldots, {\mathfrak a}_L), (f_1, \ldots, f_L))$ 
by iterated quantum Hamiltonian reduction. This algebra is conjectured to be stable equivalent to the $\CW$-algebra $\CW^k({\mathfrak g}, f)$ where $f=f_1+\cdots + f_L$. Vertex algebras are called stable equivalent if  they are isomorphic up to tensoring with free field algebras. 
Section \ref{iteratedV} extends the iterated construction to Verma modules and section 
\ref{iteratedglN} presents in detail a sequence of iterated $\CW$-algebras of 
${\mathfrak gl}_N$ associated to an ordered partition 
\[
N = m_1 s_1 + \cdots + m_L s_L 
\]
with $0 < s_1 < \dots < s_L$ and $m_1, \ldots, m_L >0$. In this case the subalgebras 
$({\mathfrak a}_1, \ldots, {\mathfrak a}_L)$ are the natural subalgebras
 $({\mathfrak g}_{N_1}, 
\ldots, {\mathfrak g}_{N_L})$, $N_i = m_i s_i$, and $f_i$ is a principal nilpotent 
element in ${\mathfrak g}_{N_i}$. The main conjecture in this paper then reads 
\begin{conj}\label{mainconj} ({\it Conjecture \ref{conjectureB}})
There is an isomorphism of vertex algebras
\be\label{eq:introE}
\cW_L^k(\mgg\gl_N, (\mgg\gl_{N_1}, \dots, \mgg\gl_{N_L}), (f_1, \dots, f_L))\cong \cW^k(\mgg\gl_N, f)\otimes  \bigotimes\limits_{1\leq i < j \leq L}  
\beta\gamma^{m_im_j(s_j-s_i)}
\ee
where each $\beta\gamma^{m_im_j(s_j-s_i)}$ factor denotes the vertex algebra of a 
pair of symplectic bosons. 
\end{conj} 
\noindent 
This conjecture is verified for characters in identity \eqref{eq:mainformulaA}, which reads
\be\label{eq:introF} 
\begin{split}
 \text{ch}[\cW_L^k&(\mgg\gl_N, (\mgg\gl_{N_1}, \dots, \mgg\gl_{N_L}), (f_1, \dots, f_L))] 
 = \ch[\cW^k(\mgg\gl_N, f)] \prod_{1\leq i < j \leq L}  
\text{ch}[ \beta\gamma^{m_im_j (s_j-s_i)}],
\end{split}
\ee
The presentation of iterated $\CW$-algebras is then concluded in Section \ref{refcharsect}, which presents the conformal weight filtration \cite{Li}  providing the resulting refined characters. In particular, it is observed that identity \eqref{eq:introF} holds for refined characters as well. 

Finally, Section \ref{comparison}, identifies the partition function \eqref{eq:introC} 
with the refined character of the universal Verma module of a vertex subalgebra 
of the tensor product 
\[
\cW^k(\mgg\gl_N, f)\otimes  \bigotimes\limits_{1\leq i < j \leq L}  
\beta\gamma^{m_im_j(s_j-s_i)}
\]
This subalgebra is obtained as a truncation of each $\beta\gamma$-system 
in the right hand side to the subalgebra $B(m_i, m_j, s_j-s_i)$ generated by the 
 fields $\beta$. 
A detailed discussion is provided in Example \ref{betagammaex}. Using this truncation, 
the following identity is proven in Proposition \ref{WZprop} 
\be\label{eq:introG} 
Z_{\bf r} = {\rm ch}_{\rm ref}^V[\cW^k(\mgg\gl_N, f)]\prod_{1\leq i<j\leq L}\
{\rm ch}_{\rm ref}[B(m_i,m_j,s_j-s_i)]
\ee
where 
\[
(r_0, \ldots, r_{\ell-1}) = \big(\, \underbrace{s_L, \ldots, s_L}_{m_L}, \ldots , 
\underbrace{s_1, \ldots, s_1}_{m_1}\, \big), \qquad 0< s_1 < \cdots s_L, \qquad 
m_i >0, \ 1\leq i \leq L.
\]
Using Conjecture \ref{mainconj}, one concludes that 
\be\label{eq:introH} 
Z_{\bf r} = {\rm ch}_{\rm ref}^V[\cW_{L,\beta}^k(\mgg\gl_N, (\mgg\gl_{N_1}, \dots, \mgg\gl_{N_L}), (f_1, \dots, f_L))] 
\ee
where 
$\cW_{L,\beta}^k(\mgg\gl_N, (\mgg\gl_{N_1}, \dots, \mgg\gl_{N_L}), (f_1, \dots, f_L))$ is the inverse image of 
the subalgebra 
\[
\cW^k(\mgg\gl_N, f)\otimes  \bigotimes\limits_{1\leq i < j \leq L}  
B(m_i, m_j, s_j-s_i) 
\]
via isomorphism \eqref{eq:introE}. Since the odd degree cohomology is trivial, this leads to an implicit identification of the cohomology of the affine Laumon space of type ${\bf r}$ to the universal 
Verma module of the vertex algebra $\cW_{L,\beta}^k(\mgg\gl_N, (\mgg\gl_{N_1}, \dots, \mgg\gl_{N_L}), (f_1, \dots, f_L))$. The main open question is whether this module structure is related to the shifted affine Yangian action constructed in \cite{M_slices,AGT_parabolic}. This will be left for future work. 

{\it Acknowledgements}. We are very grateful to Sujay Ashok for collaboration and sharing his insights with us at an early 
stage of this project. We would also like to thank Andrei Negut, Francesco Sala, Olivier Schiffmann,  Miroslav Rap\v{c}\'ak and 
Abay Zhakenov for very helpful discussions and correspondence. The work of T.C. was partially supported by NSERC Grant Number RES0048511. The work of D.-E.D. was partially supported by the NSF grant DMS-1802410.

\section{Affine Laumon spaces}\label{ALspaces} 

This section is a review of affine Laumon space following \cite{QZ}. The presentation will emphasize their construction in terms of representations of the chain-saw quiver.  

\subsection{ADHM orbifolds}\label{ADHMorb} 
Given two vector spaces $V$ and $W$, the moduli space 
of framed representations of the ADHM quiver is a GIT quotient $\CA//GL(V)$ where
\[
\CA=\{\, (A_1, A_2, I, J) \in {\rm End}(V)^{\times 2}\times {\rm Hom}(W,V)\times {\rm Hom}(W,V)\, |\, 
[A_1, A_2]+IJ=0\}.
\]
The general linear group $GL(V)$ acts
naturally on $\CA$, 
\[
g\times (A_1,A_2,I,J) \mapsto (gA_1g^{-1}, gA_2a^{-1}, gI, Jg^{-1}). 
\]
The stable locus $\CA^s\subset \CA$ is defined by 
ruling out  $A_i$-invariant subspaces $V'\subsetneq V$, $1\leq i\leq 2$, so that ${\rm Im}(I)\subset V'$. The resulting moduli space is a smooth quasi-projective variety of dimension $2rn$, 
where $r={\rm dim}(W)$ and $n= {\rm dim}(V)$. In the following we will fix a basis $(w_\alpha)$ in  $W$, $1\leq \alpha \leq r$, and 
write $(I_\alpha)$, $(J_\alpha)$ for the associated components of $I,J$. 

Any ordered partition ${\bf r}=(r_0, \ldots, r_{\ell-1})$, where
$r_0+\cdots + r_{\ell-1}=r$, $r_a\geq 0$, $\ell\geq 2$, determines a
cyclic group action $\bfmu_\ell \times \CA \to \CA$ as follows.
Let $\bfmu_\ell \to {\rm End}(W)$, be the linear action which maps the 
canonical generator $\omega= e^{2i\pi/\ell}$ to the diagonal matrix
\be\label{eq:murepA}
\omega \mapsto M = \rm{diag}(\ \underbrace{1, \ldots 1}_{r_0}, \underbrace{\omega, \ldots, \omega}_{r_1}, \ldots, \underbrace{\omega^{\ell-1}, \ldots, \omega^{\ell-1}}_{r_{\ell-1}}\ ).
\ee
Then one defines the linear action $\bfmu_\ell \times \CA \to \CA$, 
\[
\omega \times (A_1,A_2,I,J) \mapsto (A_1, \omega A_2, I M^{-1}, \omega M J). 
\]
This action commutes with the $GL(V)$-action on $\CA$ and preserves the stable locus, hence it descends to a $\bfmu_\ell$-action on the moduli space. 

The components of the $\bfmu_\ell$-fixed locus in the moduli space are indexed by 
ordered partitions ${\bf n}=(n_0,\cdots, n_{\ell-1})$, 
$n_0+\cdots+n_{\ell-1}=n$, 
$n_a\geq 0$, $0\leq a\leq \ell-1$. 
The fixed component $\CM({\bf r},{\bf n})$ corresponding to an ordered partition ${\bf n}$ consists of isomorphism classes of ADHM data $(A_1, A_2, I, J)$ satisfying 
\be\label{eq:fixedA}
A_1 = U_{\bf n} A_1 U_{\bf n}^{-1}, \qquad \omega A_2 =  U_{\bf n}A_2 U_{\bf n}^{-1}, \qquad IM^{-1} = U_{\bf n}^{-1} I, 
\qquad \omega M J = JU_{\bf n}
\ee
where 
\[
U_{\bf n} = \rm{diag}(\ \underbrace{1, \ldots 1}_{n_0}, \underbrace{\omega, \ldots, \omega}_{n_1}, \ldots, \underbrace{\omega^{\ell-1}, \ldots, \omega^{\ell-1}}_{n_{\ell-1}}\ ).
\]
As explained in \cite[Section. 2.3]{QZ}, each  component $\CM({\bf r},{\bf n})$ is identified to a moduli space of stable representations of the chain-saw quiver. Moreover, it is a smooth quasi-projective variety.

\subsection{Torus fixed points}\label{fpsection}  
The torus ${\bf T}= \IC^{\times}\times \IC^\times \times {\IC}^{\times r}$ acts on the moduli space of ADHM data by 
\[
(t_1, t_2, \zeta_1, \ldots, \zeta_r) \times (A_1, A_2,I,J) \mapsto 
(t_1A_1, t_2A_2, \zeta_\alpha^{-1}I_\alpha, t_1t_2 \zeta_\alpha J), \qquad 1\leq\alpha\leq r.
\]
The ${\bf T}$-fixed locus consists of a finite set of points in one-to-one correspondence to 
$r$-uples $\umu=(\mu_1, \ldots, \mu_r)$ of Young diagrams so that $\sum_{\alpha=1}^r |\mu_\alpha|=n$. The ${\bf T}$-fixed locus is moreover contained in the $\bfmu_\ell$-fixed 
locus.  

In order to fix conventions, note that a partition $n=n_1+\cdots+ n_k$, $n_1\geq \cdots\geq n_k$ will be identified with the set 
\[
\bigcup_{j=1}^k  \{(i,j)\in \IR^2\,|\, i\in \IZ,\ 1\leq i \leq n_j \} \subset \IR^2
\]
This set consists of $k$ rows of $n_j$ boxes each.  Note that our conventions for partitions are different from those 
used in \cite{hilblect,instcountA,instcountB,Lect_inst}.  In loc. cit. the partition associated to a set as above is determined by columns rather than 
rows. 

As explained next, the set of ${\bf T}$-fixed points contained 
in $\CM({\bf r},{\bf n})$ is naturally 
 identified with a set 
of $r$-uples $\umu$ of colored Young diagrams. 
%where each box $(i,j)\in \mu_\alpha$, $1\leq \alpha \leq r$ ,
%is assigned a color $0\leq c_\alpha(i,j) \leq \ell-1$. 
A Young diagram $\mu$ will be said to be $a$-colored, $0\leq a \leq \ell-1$, 
if and only if each box $(i,j)\in \mu$ is marked by $c_a(i,j)\in \{0, \ldots, \ell-1\}$ 
which is uniquely determined by 
\be\label{eq:colorfct}
c_a(i,j) \equiv  a -j+1 \ {\rm mod}\ \ell.
\ee
For a given ${\bf T}$-fixed point $\umu$, each partition $\mu_\alpha$, $1\leq \alpha \leq r$, will be $a(\alpha)$-colored, where $0\leq a(\alpha)\leq \ell-1$ is 
 is determined from the fixed point conditions \eqref{eq:fixedA} as explained below. 
 
Note that the ${\bfmu}_\ell$-module structure on $W$ determined by the action 
\eqref{eq:murepA} is 
\be\label{eq:murepB} 
W = \underbrace{{\bf 1}\oplus \cdots {\bf 1}}_{r_0} \oplus 
\underbrace{\Omega \oplus \cdots \Omega}_{r_1}\oplus \cdots \oplus 
\underbrace{\Omega^{\ell-1}\oplus \cdots\Omega^{\ell-1}}_{r_{\ell-1}}~,
\ee
where $\Omega$ is the canonical one dimensional representation of $\bfmu_\ell$. 
Hence 
\[
W = \bigoplus_{\alpha=1}^r \Omega^{a(\alpha)} 
\]
where $0\leq a(\alpha)\leq \ell-1$ is uniquely determined by the conditions 
\[
r_{0} +\cdots + r_{a(\alpha)-1} +1 \leq \alpha \leq r_0+\cdots + r_{a(\alpha)}
\]
for each $1\leq \alpha \leq r$. By convention, $r_{-1}=0$. 
 
Then a fixed point $\umu$ belongs to the component $\CM({\bf r},{\bf n})$ if and only if 
\be\label{eq:fixedB}
\sharp\{ (i,j) \in \mu_\alpha,\ 1\leq \alpha \leq r\, |\, c_{a(\alpha)}(i,j) = b\, \} = n_b
\ee
for each $0\leq b \leq \ell-1$. 

Let $(T_1,T_2, Z_\alpha)$, $1\leq \alpha\leq r$, be the canonical generators of the 
representation ring $R_{\bf T}$. For any $r$-uple $\umu$ let 
$T_\umu$ denote the
${\bf T}$-equivariant tangent space to the moduli space at the fixed point indexed by $\umu$. 
As shown in  \cite[Theorem 2.11]{instcountA}, the following identity holds in $R_{\bf T}$
\be\label{eq:tangentA}
T_{\umu} = \sum_{1\leq \alpha, \beta \leq r} \, Z_\beta Z_\alpha^{-1} \bigg( 
\sum_{(i,j) \in \mu_\alpha} 
 T_1^{-\mu_{\beta,j}+i} T_2^{\mu'_{\alpha,i}-j+1} + 
\sum_{(i,j) \in \mu_\beta} T_1^{\mu_{\alpha,j}-i +1} T_2^{-\mu'_{\beta,i}+j}\bigg),
\ee
where $\mu_j$ denotes the length of the $j$-th row of the partition $\mu$ while 
$\mu'_i$ denotes the height of the $i$-th column of $\mu$. 

Suppose the fixed point $\umu$ belongs to $\CM({\bf r},{\bf n})$ for some ordered partition 
${\bf n}$. Then $T_\umu$ has a natural ${\bf T}\times \bfmu_\ell$-module structure 
and the  ${\bf T}$-equivariant 
tangent space $T_\umu\CM_{\bf n}$ is given by the $\bfmu_\ell$-fixed part 
$T_\umu^{\bfmu_\ell}$. In complete analogy to \cite[Theorem 2.11]{instcountA}, the 
${\bf T}\times \bfmu_\ell$-module structure of $T_\umu$ 
is given by: 

\begin{lemm}\label{tanglemmA} 
For any $1\leq \alpha, \beta \leq r$ let $T_{\alpha, \beta}\in R_{{\bf T}\times \bfmu_\ell}$
be defined by 
\be\label{eq:tangentB}
T_{\alpha, \beta} = \Omega^{a(\beta)-a(\alpha)}Z_\beta Z_\alpha^{-1} \bigg( 
\sum_{(i,j) \in \mu_\alpha} 
 T_1^{-\mu_{\beta,j}+i} (\Omega T_2)^{\mu'_{\alpha,i}-j+1} + 
\sum_{(i,j) \in \mu_\beta}
T_1^{\mu_{\alpha,j}-i +1} (\Omega T_2)^{-\mu'_{\beta,i}+j}\bigg)~,
\ee 
where $\Omega$ is the canonical one dimensional representation of $\bfmu_\ell$.
Then 
\be\label{eq:tangentC} 
T_\umu = \sum_{1\leq \alpha, \beta \leq r} T_{\alpha, \beta}.
\ee
in the representation ring $ R_{{\bf T}\times \bfmu_\ell}$.
\end{lemm} 

\section{The generating function of Poincar\'e polynomials}\label{PPsection}

The goal of this section is to find a closed expression for the generating function of the 
Poincar\'e polynomials of the affine Laumon spaces $\CM({\bf r}, {\bf n})$ for 
fixed invariants ${\bf r}=(r_0, \ldots, r_{\ell-1})$. The computation 
will employ Atiyah-Bott-Morse localization by analogy to \cite[Theorem 3.8]{Lect_inst}. 

Let $\lambda: \IC^\times \to {\bf T}$ be a one parameter group of the form 
\[
t\mapsto (t^{M_1}, t^{M_2}, t_\alpha^{N_\alpha}) 
\]
where 
\[
M_2 >> N_1 > \cdots > N_r >> M_1 >0. 
\]
As shown in the proof of \cite[Theorem 3.7]{Lect_inst}, the resulting $\IC^\times$-action 
on the ADHM moduli space is compact. The same holds for the induced action on $\CM({\bf r},{\bf n})$ since the latter is a closed subscheme of the ADHM moduli space. Therefore the first step is the computation of the Morse index for an arbitrary fixed point $\umu\in \CM({\bf r}, {\bf n})$. 

\subsection{The Morse index}\label{morse}
Let $\umu=(\mu_1, \ldots, \mu_r)$ be a ${\bf T}$-fixed point in $\CM({\bf r},{\bf n})$. Recall that any integer $1\leq \beta \leq r$ determines a unique $0\leq a(\beta) \leq \ell-1$ so that 
\[
r_0 + \cdots + r_{a(\beta) -1}+1 \leq \beta \leq r_1+ \cdots + r_{a(\beta)-1} + r_{a(\beta)}. 
\]
As explained in Section \ref{fpsection}, this
 determines the coloring function for the partition $\mu_\beta$ to be 
\[
c_{a(\beta)}(i,j) \equiv a(\beta) -j +1 \ (\ell).
\]
For any $0\leq c\leq \ell-1$ and any $1\leq \beta \leq r$ let $n_c(\mu_\beta)$ be the 
total number of $c$-marked boxes in $\mu_\beta$. Let also 
${\rm col}(\mu_\beta)$ denote the total number of columns of $\beta$. 
Then the main result of this section is the following.

\begin{lemm}\label{morselemm}
The Morse index associated to the fixed point $\umu \in \CM({\bf r},{\bf n})^{\bf T}$ is given by 
\be\label{eq:mindexA}
w(\umu) = \sum_{\beta=1}^r w(\mu_\beta)
\ee
where 
\be\label{eq:mindexB} 
w(\mu_\beta) = \sum_{c=0}^{\ell-1} n_c(\mu_\beta) r_{c} - {\rm col}(\mu_\beta)( r_0+\cdots + r_{a(\beta)}-\beta+1).
\ee
\end{lemm} 

{\it Proof}.
 The tangent space to 
$\CM({\bf r},{\bf n})$ at the fixed point $\umu$ is the $\bfmu_\ell$-fixed part of 
the right hand side of equation \eqref{eq:tangentC},
\[
T_{\umu}^{\bfmu_\ell} = \sum_{1\leq \alpha, \beta \leq r} T^{\bfmu_\ell}_{\alpha, \beta}~.
\]
Using equation \eqref{eq:tangentB}, one obtains 
\be\label{eq:tangentD}
\bal 
T^{\bfmu_\ell}_{\alpha, \beta} = &\  
  \sum_{\substack{(i,j) \in \mu_\alpha \\ 
\mu'_{\alpha, i} - j +1 +a(\beta)-a(\alpha) \equiv 0\, (\ell)}}\ 
T_1^{-\mu_{\beta,j}+i} T_2^{\mu'_{\alpha,i}-j+1} \\
& +  \sum_{\substack{(i,j) \in \mu_\beta\\ 
-\mu'_{\beta,i}+j + a(\beta)-a(\alpha) \equiv 0\, (\ell)}}\
T_1^{\mu_{\alpha,j}-i +1} T_2^{-\mu'_{\beta,i}+j}
\eal
\ee
As in the proof of \cite[Theorem 3.8]{Lect_inst}, the Morse index of the fixed point 
$\umu$ is obtained by counting the number of monomials in the expression of 
$T_{\umu}^{\bfmu_\ell}$ satisfying the following conditions:
\begin{itemize} 
\item[$(i)$] If $\alpha \geq \beta$ then the weight of $T_2$ is negative.
\item[$(ii)$] If $\alpha < \beta$ then the weight of $T_2$ is nonpositive. 
\end{itemize} 
Note that all such monomials are contained in the second term in the right hand side of equation \eqref{eq:tangentD}. All monomials in the first term have positive $T_2$-degree.
Equation \eqref{eq:tangentD} shows that 
$w(\umu)$ has a decomposition 
\[
w(\umu)= \sum_{\beta=1}^r \sum_{\alpha=1}^r w_{\beta,\alpha}(\umu) 
\]
where 
\[
w_{\beta, \alpha} (\umu) = \left\{\begin{array}{ll} 
\sharp \{ (i,j)\in \beta \,|\, \mu'_i - j \equiv a(\beta)-a(\alpha)\ (\ell)\}, 
& {\rm for}\ \alpha <\beta,\\
& \\ 
\sharp \{ (i,j)\in \beta \,|\, j \leq \mu_i'-1, \
\mu'_i - j \equiv a(\beta)-a(\alpha)\ (\ell)\}, & 
{\rm for}\ \alpha \geq \beta.
\end{array} \right. 
\]
As shown in Lemma \ref{morselemmA}, equation \eqref{eq:fixedcountAB}, 
\be\label{eq:morsecountA}
\sharp \{ (i,j)\in \beta \,|\, \mu'_i - j \equiv a(\beta)-a(\alpha)\ (\ell)\} =
n_{a(\alpha)}(\mu_\beta). 
\ee
Furthermore, Lemma \ref{morselemmB}, equation \eqref{eq:fixedcountD}, shows that 
\be\label{eq:morsecountB} 
\sharp \{ (i,j)\in \beta \,|\, j \leq \mu_i'-1, \
\mu'_i - j \equiv a(\beta)-a(\alpha)\ (\ell)\} = \left\{\begin{array}{ll} 
n_{a(\alpha)}(\mu_\beta),  & {\rm if}\ a(\alpha) \neq a(\beta) \\
& \\ 
n_{a(\alpha)}(\mu_\beta)-{\rm col}(\mu_\beta),  & {\rm if}\ a(\alpha) = a(\beta). \\
\end{array} \right.
\ee
Then 
\[
w(\umu) = \sum_{\beta=1}^r \left(w_{<}(\mu_\beta) + w_{0}(\mu_\beta) 
+ w_{>}(\mu_\beta)\right),
\]
where 
\[
\bal
w_{<}(\mu_\beta) &\ =  \sum_{\alpha=1}^{r_0+\cdots + r_{a(\beta)-1}} 
n_{a(\alpha)}(\mu_\beta) \\
&\ = \sum_{c=0}^{a(\beta)-1} n_c(\mu_\beta)r_c~, 
\eal
\]
\[
\bal
w_{0}(\mu_\beta) &\ = \sum_{\alpha=r_0+\cdots + r_{a(\beta)-1}+1}^{\beta-1}
n_{a(\alpha)}(\mu_\beta) + 
\sum_{\alpha = \beta}^{ r_0+\cdots + r_{a(\beta)} }
\left(n_{a(\alpha)}(\mu_\beta) -{\rm col}(\mu_\beta)\right)\\
&\ = 
n_{a(\beta)}(\mu_\beta)r_{a(\beta)} - {\rm col}(\mu_\beta)( r_0+\cdots + r_{a(\beta)}-\beta+1)~,
\eal
\]
\[
\bal 
w_{<}(\mu_\beta) &\ = \sum_{\alpha = r_0+\cdots + r_{a(\beta)}+1}^r 
n_{a(\alpha)}(\mu_\beta)\\
&\ = \sum_{c=a(\beta)+1}^{\ell-1} n_c(\mu_\beta) r_{c}~.
\eal
\]
Then 
\[
\bal 
w(\mu_\beta) &\ = w_{<}(\mu_\beta) + w_0(\mu_\beta) + w_{>}(\mu_\beta) \\
&\ = \sum_{c=0}^{\ell-1} n_c(\mu_\beta) r_{c} - {\rm col}(\mu_\beta)( r_0+\cdots + r_{a(\beta)}-\beta+1).
\eal
\]

\hfill $\Box$

\subsection{The generating function}\label{genfctsect}
In order to compute the generating function, first note the
following partition identity 
derived in \cite[Section 3.2]{Kanno:2011fw}.
\begin{lemm}\label{eq:partidlemm} 
Let $v$ and $\{X_a\}$, $a\in \IZ$, be formal variables so that $X_a=X_{a+\ell}$ for 
any $a\in \IZ$. 
Then the following identity holds for any partition $\mu$ and any $a\in \IZ$. 
\be\label{eq:partidB} 
\sum_{\mu} v^{{\rm col}(\mu)} \prod_{(i,j)\in \mu} X_{a-j+1} = 
\prod_{k=1}^\infty {1\over 1- vz^k} \prod_{c=1}^{\ell-1} 
{1\over 1- vz^k \prod_{b=c}^{\ell-1} X^{-1}_{a-b}}~,
\ee
where 
\[
z=\prod_{c=0}^{\ell-1} X_c.
\]
\end{lemm}

{\it Proof}. The proof is straightforward once one notices that 
\[
v^{{\rm col}(\mu)}\prod_{(i,j)\in \mu} X_{a-j+1} = 
\prod_{i=1}^{{\rm col}(\mu)}\bigg( v \prod_{j=1}^{\mu'_i} 
X_{a-j+1}\bigg),
\]
where $\mu'_i$ denotes the height of the $i$-th column. This yields the identity 
 \be\label{eq:partidA} 
\sum_{\mu} v^{{\rm col}(\mu)} \prod_{(i,j)\in \mu} X_{a-j+1} = 
\prod_{k=1}^\infty {1\over 1- v\prod_{j=1}^k X_{a-j+1}}~.
\ee
Identity \eqref{eq:partidB} is obtained by rewriting the right hand side of \eqref{eq:partidA} 
using the relations $X_a= X_{a+\ell}$.

\hfill $\Box$  

Now let  
\[
Z_{\bf r} = \sum_{{\bf n}} P_y(\CM_{\bf n}) \prod_{a=0}^{\ell-1} q_a^{n_a} 
\]
be the generating function of Poincar\'e polynomials of affine Laumon spaces 
with fixed invariants $(r_0, \ldots, r_{\ell-1})$. For any $a\in \IZ$ let 
$q_a=q_c$ and $r_a=r_c$, where $c\in \{0, \ldots, \ell-1\}$ is uniquely determined by $a\equiv c$ mod 
$\ell$.
\begin{theo}\label{PPthmA} 
The following identity holds
\be\label{eq:PPfourB}
Z_{\bf r}=  
\prod_{a=0}^{\ell-1} \prod_{t=1}^{r_a} \prod_{c=1}^{\ell-1}
\prod_{k=1}^{\infty} {1\over 1- y^{-2t} z^{k} }
{1\over 
1- y^{-2t} z^{k}(\wq_{a-c})^{-1} (\wq_{a-c-1})^{-1} \cdots 
(\wq_{a-\ell+1})^{-1}}~,
\ee
where 
\be\label{eq:chvarA}
\wq_a = y^{2r_a} q_a
\ee
for any $a\in \IZ$, and 
 \[
 z=\prod_{a=0}^{\ell-1} \wq_a.
\]
 \end{theo} 

{\it Proof}. By
Atiyah-Bott-Morse localization, the generating function is given by 
\be\label{eq:PPone} 
\bal 
\sum_{(\mu_1, \ldots, \mu_r)} \prod_{\beta=1}^r y^{2w(\mu_\beta)} 
\prod_{c=0}^{\ell-1} q_c^{n_c(\mu_\beta)} = &\ 
\prod_{\beta=1}^r \sum_{\mu_\beta} y^{2w(\mu_\beta)} 
\prod_{c=0}^{\ell-1} q_c^{n_c(\mu_\beta)}.\\
\eal 
\ee
where $w(\mu_\beta)$ is given by \eqref{eq:mindexB}. 
For each $1\leq \beta\leq r$, the sum in the right hand side of \eqref{eq:PPone} is over all $a(\beta)$-colored partitions 
$\mu_\beta$. Note that 
\[
\bal 
\sum_{\mu_\beta} y^{2w(\mu_\beta)} 
\prod_{c=0}^{\ell-1} q_c^{n_c(\mu_\beta)} =&\  
\sum_{\mu_\beta} y^{-2{\rm col}(\mu_\beta)( r_0+\cdots + r_{a(\beta)}-\beta+1)} 
\prod_{c=0}^{\ell-1} \big(y^{2r_{c}} q_c\big)^{n_{c}(\mu_\beta)}\\
=&\ \sum_{\mu_\beta} y^{-2{\rm col}(\mu_\beta)( r_0+\cdots + r_{a(\beta)}-\beta+1)} 
\prod_{c=0}^{\ell-1} (\wq_c)^{n_{c}(\mu_\beta)}\\
&\ \sum_{\mu_\beta} y^{-2{\rm col}(\mu_\beta)( r_0+\cdots + r_{a(\beta)}-\beta+1)} 
\prod_{(i,j)\in \mu_\beta} \wq_{a(\beta)-j+1}\eal
\]
For any $1\leq \beta \leq r$ let 
\[
v_\beta = y^{2(r_0+\cdots + r_{a(\beta)}-\beta+1)}.
\]
Then identity \eqref{eq:partidB} yields 
\be\label{eq:PPthree} 
\bal 
& \sum_{\mu_\beta} y^{2w(\mu_\beta)} \prod_{c=0}^{\ell-1} q_c^{n_c(\mu_\beta)} =\\ & \prod_{k=1}^{\infty} {1\over 1- v_\beta^{-1} z^{k} }
\prod_{c=1}^{\ell-1} {1\over 
1- v_\beta^{-1} z^{k} (\wq_{a(\beta)-c})^{-1} (\wq_{a(\beta)-c-1})^{-1} \cdots 
(\wq_{a(\beta)-\ell+1})^{-1}}~.
\\
\eal
\ee
In conclusion, one obtains 
\be\label{eq:PPfour}
Z_{\bf r}= \prod_{\beta=1}^r \prod_{k=1}^{\infty} {1\over 1- v_\beta^{-1} z^{k} }
\prod_{c=1}^{\ell-1} {1\over 
1- v_\beta^{-1} z^{k}(\wq_{a(\beta)-c})^{-1} (\wq_{a(\beta)-c-1})^{-1} \cdots 
(\wq_{a(\beta)-\ell+1})^{-1}}
\ee
where 
 \[
 z=\prod_{a=0}^{\ell-1} \wq_a.
\]
Next note the set isomorphism 
\[
\{r_0+\cdots + r_{a(\beta)-1}+1\leq \beta' \leq r_0+\cdots+ r_{a(\beta)}\}
\xrightarrow{\sim} \{ 1, \ldots, r_{a(\beta)}\} 
\]
sending 
\[
\beta' \mapsto 
 r_0+\cdots + r_{a(\beta)}-\beta'+1.
\]
Since 
 $v_\beta = y^{2(r_0+\cdots + r_{a(\beta)}-\beta+1)}$, this yields the following equivalent expression for 
identity \eqref{eq:PPfour} 
\[
Z_{\bf r} =  
\prod_{a=0}^{\ell-1} \prod_{c=1}^{\ell-1}\prod_{t=1}^{r_a} 
\prod_{k=1}^{\infty} {1\over 1- y^{-2t} z^{k} }
{1\over 
1- y^{-2t} z^{k}(\wq_{a-c})^{-1} (\wq_{a-c-1})^{-1} \cdots 
(\wq_{a-\ell+1})^{-1}}~.
\]

\hfill $\Box$

The next goal is to express the generating function \eqref{eq:PPfour} in terms of the new variables 
\be\label{eq:newvar}
 u_c  = (\wq_{-c})^{-1} \cdots 
(\wq_{-\ell+1})^{-1}~, \qquad 1\leq c\leq \ell-1.
\ee
\begin{prop}\label{newvarprop} 
The following identity holds 
\be\label{eq:PPeightB}
 \bal 
  Z_{\bf r}=\ &  \bigg(\prod_{a=0}^{\ell-1} \prod_{t=1}^{r_a} \prod_{k=1}^\infty 
  {1\over 1-y^{-2t}z^k}  \bigg)  \\
  & 
\prod_{1\leq a< c\leq \ell}\ \bigg( 
\prod_{t=1}^{r_{\ell-c}}\prod_{k=1}^\infty 
{1\over 1- y^{-2t} z^{k} u_{a}u_{c}^{-1}}\bigg)\
\bigg( \prod_{t=1}^{r_{\ell-a}}\prod_{k=1}^\infty 
{1\over 1- y^{-2t} z^{k-1} u_{a}^{-1}u_{c}}\bigg),
\eal 
 \ee
 where the variables $u_c$, $1\leq c\leq \ell-1$, are defined in \eqref{eq:newvar} 
 and $u_\ell =1$. 
 \end{prop}

{\it Proof}. 
 The right hand side of \eqref{eq:PPfourB} is rewritten as follows 
 \be\label{eq:PPfourC}
 \bal 
 & \prod_{a=0}^{\ell-1} \prod_{t=1}^{r_a} \prod_{c=1}^{\ell-1} \prod_{k=1}^\infty {1\over 1-y^{-2t} z^k}
 {1\over 1-y^{-2t} z^k (\wq_{a-c})^{-1} (\wq_{a-c-1})^{-1} \cdots 
(\wq_{a-\ell+1})^{-1}} = \\
&  \prod_{a=0}^{\ell-1} \prod_{t=1}^{r_a} \prod_{k=1}^\infty {1\over 1-y^{-2t} z^k}\\
&\prod_{t=1}^{r_0}  \prod_{c=1}^{\ell-1} \prod_{k=1}^\infty  {1\over 1-y^{-2t} z^k (\wq_{-c})^{-1} (\wq_{-c-1})^{-1} \cdots 
(\wq_{-\ell+1})^{-1}} \\
&\prod_{c=1}^{\ell-1} \prod_{t=1}^{r_c} \prod_{k=1}^\infty  
{1\over 1-y^{-2t} z^k (\wq_{0})^{-1} (\wq_{-1})^{-1} \cdots 
(\wq_{c-\ell+1})^{-1}}\\
& \prod_{a=1}^{\ell-1}\prod_{t=1}^{r_a}
\prod_{\substack{c=1\\ c\neq a}}^{\ell-1} \prod_{k=1}^\infty 
{1\over 1-y^{-2t} z^k (\wq_{a-c})^{-1} (\wq_{a-c-1})^{-1} \cdots 
(\wq_{a-\ell+1})^{-1}} 
 \eal
 \ee
 Using the variables \eqref{eq:newvar},  the first two lines in the right hand of  \eqref{eq:PPfourC} become 
\be\label{eq:PPfourCB}
\bal 
  & \bigg(\prod_{a=0}^{\ell-1} 
  \prod_{t=1}^{r_a}\prod_{k=1}^\infty {1\over 1-y^{-2t} z^k}\bigg) \ 
\bigg( \prod_{c=1}^{\ell-1}\ \prod_{t=1}^{r_0} \prod_{k=1}^\infty 
 {1\over 1- y^{-2t}z^k u_c} \bigg)\\
 \eal 
 \ee
 In order to express the third line in terms of the variables \eqref{eq:newvar} note the 
 identity 
 \[
 \bal
 (\wq_{0}\wq_{-1} \cdots 
\wq_{c-\ell+1})^{-1}\ \prod_{b=0}^{\ell-1} \wq_{b} & = 
\wq_0^{-1}(\wq_{\ell-1}\cdots \wq_{c+1})^{-1} \prod_{b=0}^{\ell-1} \wq_{b} \\
&= \wq_1\cdots \wq_c \\
& = \wq_{c-\ell} \cdots \wq_{-\ell+1} \\
& = u_{\ell-c}^{-1} \\
\eal
\]
 Then the third line in the right hand side of \eqref{eq:PPfourC} 
 becomes 
 \be\label{eq:PPfourCC}
 \bal
 \bigg(\prod_{c=1}^{\ell-1} \prod_{t=1}^{r_{\ell-c}} \prod_{k=1}^\infty  
 {1\over 1- y^{-2t}z^{k-1} u_c^{-1}}\bigg)\\ 
 \eal 
 \ee
The remaining factors  are 
\be\label{eq:PPfourD} 
\bal 
& \prod_{a=1}^{\ell-1}\prod_{t=1}^{r_a}
\prod_{\substack{c=1\\ c\neq a}}^{\ell-1} \prod_{k=1}^\infty 
{1\over (1-y^{-2t} z^k (\wq_{a-c})^{-1} (\wq_{a-c-1})^{-1} \cdots 
(\wq_{a-\ell+1})^{-1})} ~.
\eal 
 \ee
As shown in Appendix \ref{infprodid} the following identity holds 
\be\label{eq:PPfourDB}
\bal
 & \prod_{a=1}^{\ell-1}\prod_{t=1}^{r_a}
\prod_{\substack{c=1\\ c\neq a}}^{\ell-1} 
{1\over (1-y^{-2t} z^k (\wq_{a-c})^{-1} (\wq_{a-c-1})^{-1} \cdots 
(\wq_{a-\ell+1})^{-1})} =\\
 & \prod_{1\leq a< c\leq \ell-1}\
\bigg( 
\prod_{t=1}^{r_{\ell-c}}
{1\over 1- y^{-2t} z^{k} u_{a}u_{c}^{-1}}\bigg)\
\bigg( \prod_{t=1}^{r_{\ell-a}}{1\over 1- y^{-2t} z^{k-1} u_{a}^{-1}u_{c}}\bigg). \\
\eal
\ee
Then, using equations \eqref{eq:PPfourCB} and \eqref{eq:PPfourCC}, 
the right hand side of \eqref{eq:PPfourC} is finally equal to 
 \be\label{eq:PPeight}
 \bal 
  & \bigg(\prod_{a=0}^{\ell-1} \prod_{t=1}^{r_a} \prod_{k=1}^\infty {1\over 1-y^{-2t}z^k}  \bigg) \
  \bigg( \prod_{a=1}^{\ell-1}\ \prod_{t=1}^{r_0} \prod_{k=1}^\infty 
 {1\over 1- y^{-2t}z^k u_a} \bigg)\\
 & 
 \bigg(\prod_{a=1}^{\ell-1} \prod_{t=1}^{r_{\ell-a}} \prod_{k=1}^\infty  
 {1\over 1- y^{-2t}z^{k-1} u_a^{-1}}\bigg)\\   
 & \prod_{1\leq a< c\leq \ell-1}\
\bigg( 
\prod_{t=1}^{r_{\ell-c}}\prod_{k=1}^\infty
{1\over 1- y^{-2t} z^{k} u_{a}u_{c}^{-1}}\bigg)\
\bigg( \prod_{t=1}^{r_{\ell-a}}\prod_{k=1}^\infty
{1\over 1- y^{-2t} z^{k-1} u_{a}^{-1}u_{c}}\bigg). \\
\eal
\ee 
 Setting $u_\ell=1$, this can be further rewritten as 
 \[
 \bal 
  & \bigg(\prod_{a=0}^{\ell-1} \prod_{t=1}^{r_a} \prod_{k=1}^\infty {1\over 1-y^{-2t}z^k}  \bigg)  
\prod_{1\leq a< c\leq \ell}\
\bigg( 
\prod_{t=1}^{r_{\ell-c}}\prod_{k=1}^\infty
{1\over 1- y^{-2t} z^{k} u_{a}u_{c}^{-1}}\bigg)\
\bigg( \prod_{t=1}^{r_{\ell-a}}\prod_{k=1}^\infty
{1\over 1- y^{-2t} z^{k-1} u_{a}^{-1}u_{c}}\bigg). 
\eal 
 \]
 as claimed in equation \eqref{eq:PPeightB}. 
 
 \hfill $\Box$ 
 
\section{$\cW$-algebras}\label{section:walgebras}

$\cW$-algebras are obtained by affine vertex algebras via quantum Hamiltonian reduction.  In this section we will study an iterated version. 

\subsection{Affine Lie algebras}\label{section:affine}

This section is a brief reminder of universal Verma modules for affine Lie algebras. 
Let $\mathfrak g$ a Lie algebra, $\mathfrak h$ a Cartan subalgebra of rank $r+1$ and $\alpha_0, \alpha_1, \dots, \alpha_r$ a set of positive simple roots.
Let 
\[
\mathfrak g = \mathfrak g_+ \oplus \mathfrak h \oplus \mathfrak g_-
\]
the usual triangular decomposition where $\mathfrak g_+$ is the sum of all positive root spaces and similarly $\mathfrak g_-$ is the sum of all negative ones. 
Let $R=\mathbb C[x_1, \dots, x_r, x_0]$ and $K=\mathbb C(x_1, \dots, x_r, x_0)$. The universal Verma module of $\mathfrak g$ is a usual Verma module over $R$ and so after base change over $K$. This means, let $h_0, \dots, h_n$ be a basis of $\mathfrak h$, e.g. the one corresponding to the simple roots $\alpha_0, \dots, \alpha_r$. 
Then one considers the one dimensional module $\mathbb C_x$ over $\mathfrak h$ on which $h_i$ acts by multiplication by $x_i$.  This extends to a 
$\mathfrak g_{\geq 0} =\mathfrak h \oplus \mathfrak g_+$-module by letting $\mathfrak g_+$ act trivially. The universal Verma module is then
\[
V := U(\mathfrak g) \otimes_{ U(\mathfrak g_{\geq 0}) } \mathbb C_x.
\]
In particular as a $\mathfrak h$-graded vector space $M$ is isomorphic to $\mathbb C[B_-]$ for $B_-$ a homogeneous basis (every basis vector is an $\mathfrak h$-eigenvector) of $\mathfrak g_-$.

Let us now specialize to $\mathfrak g$ being the affinization of a finite dimensional reductive Lie algebra $\mathfrak g^0$. 
Let $\alpha_1, \dots, \alpha_r$ be a set of positive simple roots for the finite dimensional subalgebra, let $\theta$ be the longest root and let $\delta$ be the imaginary root of $\mathfrak g$. Then $\alpha_0 = \delta-\theta$. Let $\Delta_\pm^0, \Delta^0$ denote the sets of positive, negative roots and roots of $\mathfrak g^0$ . 
Then the positive root space is 
\[
\Delta_+ = \Delta_+^0 \cup \{ \alpha+n\delta, n\delta | n\in \mathbb Z_{>0},  \alpha \in \Delta^0\}
\]
and similarly 
\[
\Delta_+ = \Delta_-^0 \cup \{ \alpha+n\delta, n\delta | n\in \mathbb Z_{<0},  \alpha \in \Delta^0\}
\]
The dimension of the root space for $\alpha + n\delta$ is one and the one of $n\delta$ is the rank of $\mathfrak g^0$. 
It follows that the formal character of the universal Verma module is given by the Weyl denominator, that is
\begin{equation}
\begin{split}
\text{ch}[V] &= \frac{\text{ch}[\mathbb C_x]}{D}   \\
D &=  \prod_{n=1}^\infty (1- e^{-n\delta})^{r} \prod_{n=1}^\infty \prod_{\alpha \in \Delta_+^0} (1-e^{-\alpha} e^{(1-n)\delta}) (1- e^{\alpha} e^{-n\delta}). 
\end{split} 
\end{equation}
This is an equation as formal power series with the short hand notation
\[
\frac{1}{1-y} = \sum_{n=0}^\infty y^n.
\]
Here 
\[
\text{ch}[\mathbb C_x] = \prod_{i=0}^r e^{x_i \omega_i}
\]
with $\omega_i$ dual to our basis of the Cartan subalgebra, that is $\omega_i(h_j) = \delta_{i, j}$.
The term universal is due to that for a ring homomorphism $f: R\rightarrow R$ (or $\mathbb C$ if desired) one can compose the action of $\mathfrak g$ with $f$ to get a module with character
\[
\text{ch}[V_f] = \prod_{i=0}^r e^{f(x_i) \omega_i}  \frac{1}{D}.
\]
in particular any ordinary Verma module over $\mathbb C$ can be realized in this way. 
If $f$ is an automorphism, then we call $M_f$ a universal Verma module as well. 
An example is that $f$ maps each $x_i$ to $x_i + \lambda_i$ for certain constants $\lambda_i$. 

We note that any ordinary Verma module is simple for generic highest-weight. Moreover it is also projective and injective in the category of weight Verma modules, that is semisimple action of the Cartan subalgebra.

\subsubsection{$\mathfrak g = \widehat{\mathfrak{gl}}_N$}\label{vermagln}

We identify $\mathfrak{gl}_N$ with its defining matrix realization. Let $E_{i, j}$ be the matrix with entry one at position $(i, j)$ and zero otherwise. 
Set $h_i = E_{i, i}$ and define $\alpha_i$ by $\alpha_i(h_j) = \delta_{i, j}$. The positive simple roots are then $\alpha_i -\alpha_{i+1}$ for $i = 1, \dots, N-1$. 
The affinization has then basis $E_{i, j}^n$ for $1\leq i, j, \leq N$ and integer $n$. The elements corresponding to $\mathfrak{g}_+$ are the 
$E_{i, j}^n$ for $1\leq i, j, \leq N$ and either positive integer $n$ or $n=0$ and $i<j$.  Similarly the elements corresponding to $\mathfrak{g}_-$ are the 
$E_{i, j}^n$ for $1\leq i, j, \leq N$ and either negative integer $n$ or $n=0$ and $i>j$.
Let $\mathbb C_x$ be as above, that is $h_i$ acts by multiplication with $x_i$ and any element in $\mathfrak{g}_+$ acts trivially. 
Let $M$ be the induced universal Verma module. Since $\omega_i = \alpha_i$ 
and since
\[
\Delta_+^0 = \{ \alpha_i - \alpha_j| 1 \leq  1 < j  \leq r \}
\]
we have
\begin{equation}
\begin{split}
\text{ch}[M] &= \frac{\text{ch}[\mathbb C_x]}{D}   \\
\text{ch}[\mathbb C_x] &= \prod_{i=0}^N e^{x_i \alpha_i}\\
D &=  \prod_{n=1}^\infty (1- e^{-n\delta})^{N} \prod_{n=1}^\infty \prod_{1 \leq i < j \leq r} (1-e^{\alpha_j - \alpha_i} e^{(1-n)\delta}) (1- e^{\alpha_i-\alpha_j} e^{-n\delta}). 
\end{split} 
\end{equation}
This converges to a meromorphic function on 
\[
\{ |v_1| > |v_2| > \dots > |v_N|\} \cap \{ |zv_i| < |v_{i+1}| | i = 1, \dots N-1\} 
\subset  \mathbb C^N \times \mathbb H
\]
 if we set $v_i = e^{\alpha_i}$ and $z = e^{-\delta}$. The expression has a meromorphic continuation to 
\begin{equation}\label{eq:univchar}
\begin{split}
\text{ch}[M] &=  \frac{\prod\limits_{i=0}^N v_i^{x_i}}{ \prod\limits_{n=1}^\infty (1- z^n)^{N} \prod\limits_{n=1}^\infty \prod\limits_{1 \leq i < j \leq N} (1-v_j v_i^{-1} z^{n-1}) 
(1- v_i v_j^{-1} z^n)}. 
\end{split} 
\end{equation}

\subsection{Quantum Hamiltonian reduction}\label{QHsection}

We summarize $\cW$-algebra basics following \cite{KWIII}  and using the summary of  \cite[Section 3]{CL}. For this let $\mgg$ be a Lie algebra,  
\begin{equation}
( \ \ | \ \ ): \mgg \times \mgg \rightarrow \mathbb C
\end{equation}
a nondegenerate invariant symmetric bilinear form on $\mgg$ and 
 $\{q^\alpha\}_{\alpha \in S}$ be a basis of $\mgg$ indexed by the set $S$.
 The structure constants are defined by 
\[
[q^\alpha, q^\beta] = \sum_{\gamma \in S}{f^{\alpha\beta}}_\gamma q^\gamma.
\]
The affine vertex algebra $V^k(\mgg)$ of $\mgg$ at level $k \in \mathbb C$  is strongly and freely generated by $\{X^\alpha\}_{\alpha \in S}$  with the usual operator products 
\[
X^\alpha(z)X^\beta(w) \sim \frac{k(q^\alpha|q^\beta)}{(z-w)^2} + \frac{\sum_{\gamma\in S}{f^{\alpha \beta}}_\gamma X^\gamma (w) }{(z-w)}.
\]
We assume that $k$ is generic. It is more precise to denote $V^k(\mgg)$ by $V^{k( \ \ | \ \ )}(\mgg)$ to emphazise the given bilinear form and we will do so later. 
Set  $X_\alpha$ to be the field corresponding to $q_\alpha$ where $\{q_\alpha\}_{\alpha \in S}$ is the dual basis of $\mgg$ with respect to $( \ \ | \ \  )$.

Let $f$ be a nilpotent element of $\mgg$. By the Jacobson-Morozov theorem, $f$ is part of an $\gs\gl_2$ triple $\{f, x, e\} \subseteq \mgg$ satisfying the  relations $[x, e] =e, [x, f]=-f, [e, f]= 2x$. 
In particular  $\mgg$ decomposes as an $\gs\gl_2$-module
\[
\mgg = \bigoplus_{k \in  \frac{1}{2}\mathbb Z} \mgg_k, \qquad \mgg_k = \{  a \in \mgg | [x, a] = ka \}. 
\] 
Let $S_k$ be a basis of $\mgg_k$ and take $S=\bigcup_k S_k$ as basis of $\mgg$. 
Let 
us also set 
\[
\mgg_+ = \bigoplus_{k \in  \frac{1}{2}\mathbb Z_{>0}} \mgg_k, \qquad \mgg_- = \bigoplus_{k \in  \frac{1}{2}\mathbb Z_{<0}} \mgg_k
\]
with corresponding basis $S_+$ of $\mgg_+$; $\mgg_-$ is naturally identified with the dual of $\mgg_+$. 
On $\mgg_{\frac{1}{2}}$ one defines the antisymmetric bilinear form
\[
\langle a, b \rangle := ( f | [a, b] ). 
\]
Let $F(\mgg_+)$ be the vertex superalgebra associated to the vector superspace $\mgg_+ \oplus \mgg_+^*$. This vertex superalgebra is strongly generated by odd fields $\{\varphi_\alpha, \varphi^\alpha\}_{\alpha \in S_+}$. The conformal weight of $\varphi_\alpha$ is $j$ if $q^\alpha$ is in $\mgg_j$ and the one of $\varphi^\alpha$ is $1-j$.  
The operator products are
\[
\varphi_\alpha(z) \varphi^\beta(w) \sim \frac{\delta_{\alpha, \beta}}{(z-w)}, \qquad \varphi_\alpha(z) \varphi_\beta(w) \sim 0 \sim \varphi^\alpha(z) \varphi^\beta(w).
\]
Let $F(\mgg_{\frac{1}{2}})$ be the neutral vertex algebra associated to $\mgg_{\frac{1}{2}}$ with bilinear form $\langle \ \ , \ \ \rangle$. This has even strong generators $\{\Phi_\alpha\}_{\alpha \in S_{\frac{1}{2}}}$ with operator products 
\begin{equation}\label{eq:neutral}
\Phi_\alpha(z) \Phi_\beta(w) \sim \frac{\langle q^\alpha , q^\beta \rangle}{(z-w)}  \sim \frac{ ( f | [q^\alpha, q^\beta] )}{(z-w)},
\end{equation}
and fields corresponding to the dual basis with respect to $\langle \ \ , \ \ \rangle$ are denoted by $\Phi^\alpha$. 
Set 
\[
F(\mgg, f, k) := F(\mgg_+) \otimes F(\mgg_{\frac{1}{2}}).
\]
The complex is 
\[
C(\mgg, f, k) := V^k(\mgg) \otimes F(\mgg, f, k).
\]
It is  $\mathbb Z$-graded by ghost number, that is the $\varphi_\alpha$ have charge minus one, the $\varphi^\alpha$ charge one and all others charge zero. 
One further defines the odd field $d(z)$ of charge minus one by
\begin{equation}
\begin{split}
d(z) &= \sum_{\alpha \in S_+}  :X^\alpha \varphi^\alpha:  
-\frac{1}{2} \sum_{\alpha, \beta, \gamma \in S_+}  {f^{\alpha \beta}}_\gamma :\varphi_\gamma \varphi^\alpha \varphi^\beta: + 
 \sum_{\alpha \in S_+} (f | q^\alpha) \varphi^\alpha + \sum_{\alpha \in S_{\frac{1}{2}}} :\varphi^\alpha \Phi_\alpha:.
\end{split} 
\end{equation}
The zero-mode $d_0$ is a differential since $[d(z), d(w)]=0$ by \cite[Thm. 2.1]{KRW}.  The $\mathcal W$-algebra is  defined to be its homology
\[
\cW^k(\mgg, f) := H\left(C(\mgg, f, k), d_0 \right). 
\]
The homology is non-trivial only in degree zero \cite{KWIII}.

\subsection{Examples for $\mgg = \mgg\gl_N$}\label{examples}

$\cW$-algebras of $\mgg\gl_N$ are parameterized by partitions of $N$. We introduce them step by step. 

\begin{exam} {\textup{(Principal $\cW$-algebra: $N= N$)}}

Let $\mgg= \mgg\gl_N$, denote the standard representation of it by $\mathbb C^N$ and its conjugate by $\overline{\mathbb C}^N$. 
Let $f$ be the principal nilpotent element. That is, let 
$\rho : \gs\gl_2 \rightarrow \mgg\gl_N = \mgg\gl(\mathbb C^N)$ be the $N$-dimensional irreducible representation of $\gs\gl_2$. We denote this by $\mathbb C_N$.  More generally, let $\IC_d$, $d\geq 1$, denote the $d$-dimensional irreducible representation 
of $\gs\gl_2$. Then $\mgg\gl_N$ decomposes as a module for the $\gs\gl_2$-subalgebra whose basis is the triple $(e, h, f)$ as 
\[
\mgg\gl_N \cong_{\mgg\gl_N}  \mathbb C^N \otimes \overline{\mathbb C}^N \cong_{\gs\gl_2} \mathbb C_N \otimes \mathbb C_N \cong_{\gs\gl_2}  
 \bigoplus_{J=1}^{N} \mathbb C_{2J-1}.
\]
The free fermion algebra associated to the quantum Hamiltonian reduction is the free fermion algebra corresponding to the non-zero weight spaces of this decompositon, whose dimension is $N(N-1)$, i.e. it is $F^{N(N-1)}$ the vertex superalgebra of $N(N-1)$ free fermions. 
The complex for the reduction is then $V^k(\mgg\gl_N) \otimes F^{N(N-1)}$. 
\end{exam}

\begin{exam} {\textup{(Rectangular $\cW$-algebra: $N = s+ s + \dots  +  s$ )}}

Let $\mgg = \mgg\gl_N$ with $N=\ell s$. There is then an embedding of $\mgg\gl_\ell \oplus \mgg\gl_s$ in $\mgg\gl_N$, such that $\mathbb C^N \cong \mathbb C^\ell \otimes \mathbb C^s$ is the tensor product of the standard representations of $\mgg\gl_\ell$ and $\mgg\gl_s$. Hence 
\[
\mgg\gl_N \cong_{\mgg\gl_N}  \mathbb C^N \otimes \overline{\mathbb C}^N \cong (\mathbb C^\ell \otimes \overline{\mathbb C}^\ell) \otimes (\mathbb C^s \otimes \overline{\mathbb C}^s)  \cong \textup{adj}(\mgg\gl_\ell) \otimes \textup{adj}(\mgg\gl_s) 
\]
is just the tensor product of the adjoint representations of $\mgg\gl_\ell$ and $\mgg\gl_s$. Taking $f$ and $\gs\gl_2$ to be the image in $\mgg\gl_N$ of the principal nilpotent element in $\mgg\gl_s$  and its associated $\gs\gl_2$ via this embedding we get from the previous example that 
\be\label{eq:spindecompA}
\mgg\gl_N \cong \textup{adj}(\mgg\gl_\ell)  \otimes  \left( \bigoplus_{J=1}^{s} \mathbb C_{2J-1} \right)
\ee
as a module for $\mgg\gl_\ell \otimes \gs\gl_2$. In particular the complex for this reduction is $V^k(\mgg\gl_N) \otimes F^{\ell^2 s(s-1)}$.
\end{exam}

\begin{exam}\label{exam2} {\textup{($\cW$-algebra associated to $N= s_1+ s_1 + \dots  +  s_1 + s_2+ s_2 + \dots  +  s_2$)}}

Let $N = N_1 +  N_2$ and $N_1 = m_1 s_1, N_2 = m_2 s_2$, and assume 
$0<s_1 < s_2$. We then have an embedding of $\mgg\gl_{N_1} \oplus \mgg\gl_{N_2}$ in $\mgg\gl_N$ and can choose $f = f_1 + f_2$ with $f_1$ and $f_2$ the images of the rectangular nilpotent elements in $\mgg\gl_{N_1}$ and $\mgg\gl_{N_2}$ in $\mgg\gl_N$. Then as a
$\mgg\gl_{m_1} \otimes \mgg\gl_{m_2} \otimes \gs\gl_2$-module ($\mathbf 1_v$ denots the trivial $\mgg\gl_v$-module)
\begin{equation}\label{eq:spindecompB}
\begin{split}
\mgg\gl_N
\cong\  &\, \textup{adj}(\mgg\gl_{m_1})  \otimes \mathbf 1_{m_2}  \otimes  \left( \bigoplus_{J=1}^{s_1} \mathbb C_{2J-1} \right) \oplus  \mathbf 1_{m_1}  \otimes
\textup{adj}(\mgg\gl_{m_2})  \otimes  \left( \bigoplus_{J=1}^{s_2} \mathbb C_{2J-1} \right) \oplus \\
&\ \left(\mathbb C^{m_1} \otimes \overline{\mathbb C}^{m_2} \oplus  
\overline{\mathbb C}^{m_1} \otimes {\mathbb C}^{m_2}\right) \otimes \left( \bigoplus_{\substack{d = s_2-s_1 + 1\\ d = s_2-s_1+1\ (2)}}^{s_1+s_2-1} \mathbb 
C_d\right).
\end{split}
\end{equation}
The complex for this reduction is 
\[
V^k(\mgg\gl_N) \otimes F^{m_1^2s_1(s_1-1)} \otimes F^{m_2^2s_2(s_2-1)} \otimes F^{2m_1m_2 s_1(s_2-1)}\]
 if $s_1 = s_2 \mod 2$ and  
 \[
 V^k(\mgg\gl_N) \otimes F^{m_1^2 s_1(s_1-1)} \otimes F^{m_2^2 s_2(s_2-1)} \otimes F^{2m_1m_2 s_1s_2} \otimes \beta\gamma^{m_1m_2s_1}
 \]
 otherwise. Here $\beta\gamma$ denotes the vertex algebra of a pair of symplectic bosons ($\beta\gamma$-VOA, the neutral free field algebra associated to $\mgg_{\frac{1}{2}}$).
\end{exam}

\begin{exam}\label{exam} {\textup{($\cW$-algebra associated to a general partition of $N$)}}
The general case is now $N = N_1 + \dots + N_L$, $L\geq 1$,  and $N_i = m_i s_i$, $1\leq i\leq L$, so that $\mgg\gl_N$ has $\mgg\gl_{N_1} \oplus \dots \oplus \mgg\gl_{N_L}$ as subalgebras and we can take $f = f_1 + \dots +f_L$ with $f_i$ the rectangular element in the $i$-th summand. Assume that $0<s_1 < s_2 < \dots < s_{L}$, i.e. the partition is ordered. 
Then
\be\label{eq:spindecompC}
\bal
\mgg\gl_N \cong\ &\ \  \bigoplus_{i=1}^L \textup{adj}(\mgg\gl_{m_i}) \otimes   
\left( \bigoplus_{J=1}^{s_i} \mathbb C_{2J-1} \right) \oplus
\\
& \bigoplus_{1 \leq i < j \leq L}  \left(\mathbb C^{m_i} \otimes \overline{\mathbb C}^{m_j} \oplus  \overline{\mathbb C}^{m_i} \otimes {\mathbb C}^{m_j}\right) \otimes \Bigg(\bigoplus_{\substack{d= s_j-s_i + 1\\ 
d = s_j - s_i +1\ (2)}}^{s_i+s_j-1}
 \mathbb C_d\Bigg)\\
 \eal
\ee
and 
the free field algebra of the complex is
\be\label{eq:ffalgA}
F = \bigotimes_{i=1}^L F^{m_i^2 s_i(s_i-1)} \otimes \bigoplus_{1 \leq  i < j \leq L} 
\widetilde F^{2m_im_j s_i(s_j-1)}
\ee
with 
\be\label{eq:ffalgB}
\widetilde F^{2m_i m_j s_i(s_j-1)} := \begin{cases}  
F^{2m_i m_j s_i(s_j-1)} & s_i = s_j \mod 2 \\ 
F^{2m_i m_j s_is_j} \otimes \beta\gamma^{m_i m_j s_i} & s_i \neq s_j \mod 2.
\end{cases} 
\ee
\end{exam}

\subsection{Properties}\label{properties}
\begin{enumerate} 
\item 
Set 
\begin{equation}\label{eq:Ialpha}
I^\alpha := X^\alpha + \sum_{\beta, \gamma \in S_+} {f^{\alpha \beta }}_\gamma :\varphi_\gamma\varphi^\beta: + \frac{1}{2} \sum_{\beta \in S_{\frac{1}{2}}}  {f^{\beta\alpha}}_\gamma \Phi^\beta\Phi_\gamma
\end{equation}
for $\alpha \in \mgg_0$. Denote by $\mgg^f$ the centralizer of $f$ in $\mgg$, and set $\gb := \mgg^f \cap \mgg_0$. It is a Lie subalgebra of $\mgg$.  \cite[Thm 2.1]{KWIII} says in particular 
that  $d_0(I^\alpha)=0$ for $q^\alpha \in \gb$ and
\[
I^\alpha {}(z) I^\beta(w)  \sim \frac{B(q^\alpha | q^\beta) }{(z-w)^2} + \frac{ {f^{\alpha \beta}}_\gamma I^\gamma}{(z-w)} ,
\]
with the bilinear form $B$ on $\ga$ given  by 
\[
B(q^\alpha|q^\beta)  = k(q^\alpha | q^\beta) + \frac{1}{2}\left(\kappa_\mgg(q^\alpha, q^\beta) -\kappa_{\mgg_0}(q^\alpha, q^\beta)-\kappa_{\frac{1}{2}}(q^\alpha, q^\beta) \right) 
\]
with $\kappa_\mgg$ and $\kappa_{\mgg_0}$ the usual Killing form on $\mgg$ and $\mgg_0$ and 
 $\kappa_{\frac{1}{2}}$ the trace of $\mgg_0$ on $\mgg_{\frac{1}{2}}$. 
Thus  $\cW^k(\mgg, f)$ contains $V^B(\gb)$ as subalgebra for some bilinear form depending on $(\mgg, f, k)$. 
\item 
Let $\cA^k$ be some vertex algebra containing $V^k(\mgg)$ as subalgebra, then
\[
H\left(\cA^k \otimes F(\mgg, f, k), d_0\right)
\]
is a vertex algebra that for generic $k$ contains the $\cW$-algebra $\cW^k(\mgg, f)$ as subalgebra. 
If $\cA^k$ is an object in the category of ordinary  $V^k(\mgg)$-modules, usually denoted by $\KL_k(\mgg)$ or $\KL_{k( \ \ | \ \ )}(\mgg)$, then for generic level the homology is concentrated in degree zero since Theorem 6.2. of \cite{KWIII} applies. In particular the graded character of the homology is just the Euler-Poincar\'e character, that is the supercharacter of the complex.  
\item
Let $M$ be an object in $\KL_k(\mgg)$, then $M \otimes F(\mgg, f, k)$ is a 
$V^B(\gb)$-module. The zero-modes of the $I^\alpha$ for $q^\alpha \in \gb$ give rise to an action of $\gb$ on $M \otimes F(\mgg, f, k)$. 
This action commutes with the action of the Virasoro zero-mode of  $\cW^k(\mgg, f)$ by \cite[Thm 2.1]{KWIII}. 
 Since conformal weight spaces of $M \otimes F(\mgg, f, k)$ are finite-dimensional (the conformal weight of the $\Phi^\alpha$ is $1/2$) it follows that each conformal weight space of $M \otimes F(\mgg, f, k)$ is an integrable $\gb$-module and hence $M \otimes F(\mgg, f, k)$ must be an object in a completion of $\KL_B(\gb)$. 
 \item  Let $V$ be the universal Verma module as discussed in section \ref{section:affine}. As discussed in that section, we can twist the action of $R$ on $V$ by an automorphism $\varphi$ of $R$, we call the resulting module still a universal Verma module. Let $\mathcal V^k$ be a possibly infinite direct sum of universal Verma modules (all possibly twisted).
 Then  the homology
\[
H\left(\mathcal V^k \otimes F(\mgg, f, k), d_0\right)
\]
is concentrated in degree zero, since Theorem 6.2. of \cite{KWIII} applies as well. 
 \item Let $V$ be the universal Verma module, then $V \otimes F(\mgg, f, k)$ is a 
$V^B(\gb)$-module (over $R$).  As before, the zero-modes of the $I^\alpha$ for $q^\alpha \in \gb$ give rise to an action of $\gb$ on $V \otimes F(\mgg, f, k)$. 
This action commutes with the action of the Virasoro zero-mode of  $\cW^k(\mgg, f)$. 
Any conformal weight space $V_n$ of $V$ is a direct sum of universal Verma modules for $\gb$.  Conformal weight spaces  $F(\mgg, f, k)_m$ of $F(\mgg, f, k)$  are integrable $\gb$-modules and the
tensor product of an integrable module with a universal Verma module is a direct sum of universal Verma modules.
As a $\gb$-module the weight space of conformal weight $N$ is just the finite direct sum $\bigoplus_{n} V_n \otimes F(\mgg, f, k)_{N-n}$  
and hence a direct sum of universal Verma modules for $\gb$. The action of the Cartan subalgebra of $\gb$ is semisimple on $V$ and $F(\mgg, f, k)$  and so it also is on the tensor product of the two, i.e.  $V \otimes F(\mgg, f, k)$. is an object in the category $\mathcal O$.
Let $U$ be an indecomposable  $V^B(\gb)$-submodule of  $V \otimes F(\mgg, f, k)$. Since conformal weight of the latter is lower bounded the same must be true for $U$. The top level of $U$ must be a $\gb$-module and hence must contain a universal  $\gb$-Verma module as direct summand. 
It follows that $U$ must contain the corresponding universal $V^B(\gb)$-module, but universal Verma modules are simple, projective and injective in the category $\mathcal O$, see e.g. \cite{Humphreys-book}.
We conclude that $U$ coincides with the universal $V^B(\gb)$-module.
\end{enumerate}

\subsection{Characters}\label{characters}

The Virasoro field $L$ of $\cW^k(\mgg, f)$  is
\begin{equation}
\begin{split}
L_{\text{sug}} &=  \frac{1}{2(k+h^\vee)} \sum_{\alpha \in S}  :X_\alpha X^\alpha:,\\
L_{\text{ch}} &= \sum_{\alpha \in S_+} \left(-m_\alpha :\varphi^\alpha \partial \varphi_\alpha: + (1-m_\alpha) :(\partial\varphi^\alpha )\varphi_\alpha:  \right),\\
L_{\text{ne}} &= \frac{1}{2} \sum_{\alpha \in S_{\frac{1}{2}}} :(\partial \Phi^\alpha) \Phi_\alpha : ,\\
L &= L_{\text{sug}} + \partial x + 
L_{\text{ch}} + L_{\text{ne}}.
\end{split}
\end{equation}
It has central charge 
\begin{equation}\label{eq:c}
c = c(\mgg, f, k) = \frac{k\, \text{dim}\, \mgg}{k+h^\vee} -12 k(x|x) -\sum_{\alpha \in S_+} (12m_\alpha^2-12m_\alpha +2) -\frac{1}{2}\, \text{dim}\, \mgg_{\frac{1}{2}}.
\end{equation}

Let $\mathfrak h^\sharp$ be a Cartan subalgebra of $\gb$. Choose a Cartan subalgebra  $\mfh$ of $\mgg_0$ and hence of $\mgg$ that contains $\mfh^\sharp$ and $x$. Let $\Delta_j$ be the roots in $\mgg_j$ and fix a set of positive roots $\Delta_+ = \Delta_+^0 \cup \bigcup_{j > 0} \Delta_j$ with $\Delta_+^0$ a set of positive roots for $\mgg_0$. Set $\delta(\alpha) =j$ for $\alpha\in \Delta_j$ and set $\epsilon(\alpha) = 0$ for $j \in \mathbb Z$ and $\epsilon(\alpha) = \frac{1}{2}$ otherwise. 
Let $\Delta_{h.w}$  denote those roots that correspond to a highest-weight vector for the $\gs\gl_2$-triple $(f, x, e)$. 
Let $\Delta_{h.w}^\sharp$ be the subset of those roots that act non-trivially on $h^\sharp$ and 
$\Delta_{h.w.}^0$ its complement; set
 $\Delta_{h.w}^{\pm, \sharp} = \Delta_\pm \cap \Delta_{h.w.}^\sharp$.  

The character of a module is defined to be 
\begin{equation}
\text{ch}[M](h, q) = \tr_{M}\left( q^{L_0- \frac{c}{24}} e^h\right)
\end{equation}
for  $h \in \mfh^\sharp$. 
The character of the $\cW$-algebra is
\begin{equation}\label{eq:chW}
\begin{split}
\text{ch}[\cW^k(\mgg, f)](h, q) 
&= q^{-\frac{c}{24}}   \prod_{n=1}^\infty
\prod_{\alpha \in \Delta_{\text{ h.w.}}^{0}} (1-q^{n+\delta(\alpha)})^{-1}
 \prod_{\alpha \in \Delta_{\text{ h.w.}}^{+, \sharp}} (1-e^{\alpha(h)}q^{n+\delta(\alpha)})^{-1} (1-e^{-\alpha(h)}q^{n+\delta(\alpha)})^{-1},
\end{split}
\end{equation}
where the convention is always that
\[
\frac{1}{1-y} = \sum_{n=0}^\infty y^n.
\]
The character of the Verma module $V_\lambda$ that is obtained via quantum Hamiltonian reduction from the Verma module of $V^k(\mgg)$ of highest-weight 
$\lambda$ is 
\begin{equation}\label{eq:EP}
\text{ch}[V_\lambda]  = q^{-\frac{c}{24}} q^{\frac{(\lambda| \lambda + 2\rho) }{2(k+h^\vee)} -\lambda(x)}\ \prod_{n=1}^\infty \prod_{\alpha \in \Delta_{\text{ h.w.}}^{0}} (1-q^{n-\epsilon(\alpha)})^{-1}
 \prod_{\alpha \in \Delta_{\text{ h.w.}}^{+, \sharp}}(1-e^{\alpha(h)}q^{n-\epsilon(\alpha)})^{-1} (1-e^{-\alpha(h)}q^{n+\epsilon(\alpha)-1})^{-1},
\end{equation}

Consider our main example, that is example \ref{exam}.
In this case $\gb = \mgg\gl_{m_1} \oplus \dots \oplus \mgg\gl_{m_L}$. 
We first explain how to modify the character so that the $\epsilon(\alpha)$ contributions disappear. Note that $\epsilon(\alpha) = 0$ for $\alpha \in \Delta_{\text{ h.w.}}^{0}$.
We identify $\gb$ with its natural image in $\mgg\gl_{m_1 + \dots + m_L}$, so that $\mfh^\sharp$ is identified with diagonal matrices 
whose natural basis are the $E_a$, the diagonal matrices with entry one at position $a$ and all others zero. 
$\epsilon_a$ is then the map from $\mfh^\sharp$ to $\mathbb C$ that maps $E_a$ to one and all others to zero.  Set $u_a = e^{\epsilon_a(h)}$ for $h \in \mfh^\sharp$. 
Set $x_i = \sum_{k=1}^{m_i} E_{m_1 + \dots m_{i-1} + k}$. 
Set $t_j = 0$ if $s_j = s_1 \mod 2$ and $-\frac{1}{2}$ other wise. Set $y = \sum_{k =1}^K t_k x_k$. For any monomial $\mu$ of $(u_a^{\pm 1}, y^{\pm 1}, z)$, not identically $1$, 
let  
\[
(\mu)_\infty = \prod_{n=1}^\infty(1-z^{n-1}\mu).
\]
Let $I^y$ be the field corresponding to $y$ in $\cW^k(\mgg, f)$. Then $L + \partial I^y$ is another Virasoro field. The conformal weight of a generator corresponding to $\alpha$ is unchanged if $\epsilon(\alpha)=0$ while otherwise it changes by $\pm \frac{1}{2}$ for a generator corresponding  $\Delta_{h.w}^{\pm, \sharp}$. 
Under this procedure the central charge changes to 
\[
c =  c(\mgg, f, k) - 12k(y|y)
\] 
and the Verma character becomes
\begin{equation}\label{eq:EP}
\text{ch}[V_\lambda]  = q^{-\frac{c}{24}} q^{\frac{(\lambda| \lambda + 2\rho) }{2(k+h^\vee)} -\lambda(x+y)}\ \prod_{n=1}^\infty \prod_{\alpha \in \Delta_{\text{ h.w.}}^{0}} (1-q^{n})^{-1}
 \prod_{\alpha \in \Delta_{\text{ h.w.}}^{+, \sharp}}(1-e^{\alpha(h)}q^{n})^{-1} (1-e^{-\alpha(h)}q^{n-1})^{-1},
\end{equation}
More explicitly this becomes
\begin{equation}\label{eq:charA}
\text{ch}[V_\lambda]  = q^{-\frac{c}{24}} q^{\frac{(\lambda| \lambda + 2\rho) }{2(k+h^\vee)} -\lambda(x+y)}\ \prod_{i=1}^L X_i\ \prod_{1\leq i <j \leq L} X_{i,j} \ \Big\vert_{y=1}
\end{equation}
where 
\be\label{eq:charB}
\bal
X_i = & \prod_{t=1}^{s_i} { (y^{-2t}z)_\infty^{-m_i}} \prod_{m_1+\cdots + m_{i-1} +1 \leq a<c\leq 
m_1+\cdots+ m_i } \ \ (y^{-2t} z u_au_c^{-1})^{-1}_\infty \ 
(y^{-2t} u_a^{-1} u_c)^{-1}_\infty\\
\eal 
\ee
and 
\be\label{eq:charC} 
\bal
X_{i,j} = & \prod_{t=1}^{s_i} \ \ \prod_{\substack{ m_1+\cdots + m_{i-1} +1 \leq a\leq 
m_1+\cdots+ m_i  \\ m_1+\cdots + m_{j-1} +1 \leq c \leq 
m_1+\cdots+ m_j }} \ \ (y^{-2t} z u_au_c^{-1})^{-1}_\infty  
(y^{-2t} u_a^{-1} u_c)^{-1}_\infty~.\\
\eal 
\ee

\subsection{Iterated reductions}\label{iteratedW} 

Let $\mgg$ be a simple Lie algebra and $\ga_1, \dots, \ga_L$ an ordered set of Lie algebras such that
$\ga_1 \oplus \ga_2 \dots \oplus \ga_L$ is a subalgebra of $\mgg$. Let $f_1, \dots, f_L$ be nilpotent elements in $\ga_1, \dots, \ga_L$ and let $f = f_1 + f_2 + \dots + f_L$ be their sum in $\mgg$. Each $f_i$ is also considered to be an element in $\mgg$ via the embedding of $\ga_i$ in $\mgg$. 
Define Lie algebra $\gb_1, \dots, \gb_L$ by setting $\gb_1 = \mgg$ and $\gb_{i+1} := \gb_i^{\ga_i}$ in particular each $\gb_i$ has $\ga_i \oplus \dots \oplus \ga_L$ as subalgebra and in particular $f_i$ is a nilpotent element in $\gb_i$ corresponding to an $\gs\gl_2$-triple $(e_i, h_i, f_i)$ in $\gb_i$. Let 
\[
\gb_i =  \bigoplus_{n\in \frac{1}{2}\mathbb Z} \gb_{i, n}
\]
be the corresponding grading by $h_i$ eigenvalue. 
Since $f_i \in \ga_i$ as well we have that $\gb_{i+1} =  \gb_i^{\ga_i}\subset \gb_i^{f_i} \cap \gb_{i, 0}$.
We define a sequence of vertex algebras 
$\cW^k_i(\mgg, (\ga_1, \dots, \ga_L), (f_1, \dots, f_L))$ for $i = 1, \dots, h$ with the properties that this algebra contains 
an affine vertex algebra $V^{B_{i+1}}(\gb_{i+1})$ as subalgebra for some generic level $B_{i+1}$(here set $\gb_{s+1} = 0$) and such that $\cW^k_i(\mgg, (\ga_1, \dots, \ga_L), (f_1, \dots, f_L))$ is in a completion of $\KL_{B_{i+1}}(\gb_{i+1})$. Set
\[
\cW^k_1(\mgg, (\ga_1, \dots, \ga_L), (f_1, \dots, f_L)) := \cW^k(\gb_1, f_1)
\]
and 
\[
\cW^k_i(\mgg, (\ga_1, \dots, \ga_L), (f_1, \dots, f_L)) := H\left(\cW^k_{i-1}(\mgg, (\ga_1, \dots, \ga_L), (f_1, \dots, f_L))\otimes F(\gb_{i}, f_i, B_{i}), d^i_0\right).
\]
with $d_0^i$ the differential associated to the quantum Hamiltonian reduction of $V^{B_{i}}(\gb_{i})$ with respect to the nilpotent element $f_i$.
The claimed properties hold by Properties (1,\ 2,\ 3) of the last Section \ref{properties}. In particular each homology is concentrated in degree zero. 
By construction each $d_0^i$ commutes with $d_0^j$ viewed as operators on 
\[
V^k(\mgg) \otimes F(\mgg, (\ga_1, \dots, \ga_L), (f_1, \dots, f_L))
\]
with 
\[
F(\mgg, (\ga_1, \dots, \ga_L), (f_1, \dots, f_L)) := \bigotimes_{i=1}^L 
F(\gb_{i}, f_i, B_{i})
\]
and so the iterated homology coincides with the total homology, that is
\be\label{eq:totalcpxA}
\cW^k_L(\mgg, (\ga_1, \dots, \ga_L), (f_1, \dots, f_L)) = H(V^k(\mgg) \otimes F(\mgg, (\ga_1, \dots, \ga_L), (f_1, \dots, f_L)), d_0)
\ee
and $d = d^1 + d^2 + \dots + d^L$. 
It is natural to ask if there is a relation between 
$H(V^k(\mgg) \otimes F(\mgg, (\ga_1, \dots, \ga_L), (f_1, \dots, f_L)), d_0)$ and the usual quantum Hamiltonian reduction $\cW^k(\mgg, f)$. Two vertex algebra $V_1, V_2$ are called stable equivalence if there exist a free field algebras $F_1, F_2$ such that
\[
V_1 \otimes F_1 \cong V_2 \otimes F_2
\]
as vertex algebras. Note that conformal vectors are allowed to differ. 
\begin{conj}\label{conjectureA} 
$H(V^k(\mgg) \otimes F(\mgg, (\ga_1, \dots, \ga_L), (f_1, \dots, f_L)), d_0)$ and  $\cW^k(\mgg, f)$ are stable equivalent.
\end{conj}
The graded characters of these two algebras are given by the supercharacters of the complexes which agree upon multiplication with suitable free field algebra characters. Moreover it is a general theme that quantum Hamiltonian reductions corresponding to the same nilpotent element but to different complexes and differentials coincide up to some free field algebras, see \cite{Creutzig:2018ltv} and in particular \cite{Arakawa:2020oqo} for examples. These examples are motivated from concatenation of corner vertex algebras  \cite{Gaiotto:2017euk, Creutzig:2017uxh} and the special example of our iterated reduction where $f_1, \dots, f_L$ are principal nilpotent in $\ga_1, \dots, \ga_L$ is exactly the gluing of several corner vertex algebras \cite{Prochazka:2017qum}.

\subsection{Iterated Verma modules}\label{iteratedV}

We continue with the same set-up as in the previous subsection. 
Let $M$ be a module for $V^k(\mgg)$. 
Then we can set
\[
M_1(\mgg, (\ga_1, \dots, \ga_L), (f_1, \dots, f_L))  :=  H(V^k(\mgg) \otimes F(\mfg, f_1, k), d_0^1)
\]
and iteratively
\[
M_i(\mgg, (\ga_1, \dots, \ga_L), (f_1, \dots, f_L)) := H\left(M_{i-1}(\mgg, (\ga_1, \dots, \ga_L), (f_1, \dots, f_L))\otimes F(\gb_{i}, f_i, B_{i}), d^i_0\right).
\]
Then $M_L(\mgg, (\ga_1, \dots, \ga_L), (f_1, \dots, f_L))$ is automatically a $\cW^k_L(\mgg, (\ga_1, \dots, \ga_L), (f_1, \dots, f_L))$-module. 
If $M$ is the universal Verma module $V$, then all iterated homologies are concentrated in degree zero by Properties (4) and (5) of Section \ref{properties}. Write $V^W$ for $V_L(\mgg, (\ga_1, \dots, \ga_L), (f_1, \dots, f_L))$,
In particular the homology coincides with the total homology and the character is given by the Euler-Poincar\'e character of the complex,
\begin{equation}\label{EP}
\text{ch}[V^W] 
  = \text{sch}[V \otimes F'] 
\end{equation}

\subsection{$\cW_L^k(\mgg\gl_N, (\mgg\gl_{N_1}, \dots, \mgg\gl_{N_L}), 
(f_1, \dots, f_L))$: The iterated $\cW$-algebras of $\mgg\gl_N$}\label{iteratedglN}

\begin{exam}\label{exam2} {\textup{$\cW_2^k(\mgg\gl_N, (\mgg\gl_{N_1}, \mgg\gl_{N_2}), (f_1, f_2))$}}

We reconsider Example \ref{exam}. Let $N = N_1 + N_2$ and $M_1 = m_1s_1, M_2 = m_2s_2$, and assume $0<s_1 < s_2$. We then have an embedding of $\mgg\gl_{M_1} \oplus \mgg\gl_{M_2}$ in $\mgg\gl_N$ and can choose $f = f_1 + f_2$ with $f_1$ and $f_2$ the images of the rectangular nilpotent elements in $\mgg\gl_{M_1}$ and $\mgg\gl_{M_2}$ in $\mgg\gl_M$. 
Let $\gs\gl_2^1$ be the $\gs\gl_2$ corresponding to the rectangular embedding in 
$\mgg\gl_{M_1}$ and then the latter embedded in $\mgg\gl_{N}$. 
Then as a
$\mgg\gl_{m_1} \otimes \mgg\gl_{M_2} \otimes  \gs\gl^1_2$-module ($\mathbf 1_m$ denots the trivial $\mgg\gl_m$-module)
\begin{equation} 
\bal
\mgg\gl_M \cong \ &\  \textup{adj}(\mgg\gl_{m_1})  \otimes \mathbf 1_{M_2}  \otimes  \left( \bigoplus_{J=1}^{s_1} \mathbb C_{2J-1} \right) \oplus  \mathbf 1_{m_1}  \otimes
\textup{adj}(\mgg\gl_{M_2})  \otimes \mathbf 1_{\gs\gl_2^1}  \oplus \\
&\ \left(\mathbb C^{m_1} \otimes \overline{\mathbb C}^{M_2} \oplus  
\overline{\mathbb C}^{m_1} \otimes {\mathbb C}^{M_2}\right) \otimes   \mathbb C_{s_1}
\eal
\end{equation}
The complex for this reduction is 
\[
V^k(\mgg\gl_M) \otimes F^{m_1^2 s_1(s_1-1)} \otimes F^{2m_1m_2 s_2(s_1-1)}
\]
 if $s_1$ is odd and  
 \[
 V^k(\mgg\gl_M) \otimes F^{m_1^2 s_1(s_1-1)} \otimes F^{2m_1m_2 s_1s_2} \otimes \beta\gamma^{m_1m_2s_2}
 \]
 otherwise. Here $\beta\gamma$ denotes the vertex algebra of a pair of symplectic bosons ($\beta\gamma$-VOA, the neutral free field algebra associated to $\mgg_{\frac{1}{2}}$). 
The cohomology still has an affine VOA of type $\mgg\gl_{M_2}$ and the iterated reduction is the rectangular one on this affine subVOA, in particular the complex for the combined reductions is then
\[
V^k(\mgg\gl_N) \otimes F^{m_1^2 s_1(s_1-1)} \otimes F^{m_2^2 s_2(s_2-1)}  \otimes 
\begin{cases}  F^{2m_1m_2 s_2(s_1-1)} & \quad s_1 \ \text{odd} \\ 
F^{2m_1m_2 s_1s_2} \otimes \beta\gamma^{m_1m_2s_2} & \quad s_1 \ \text{even} \end{cases}
\]
\end{exam}

The general case is now $N = N_1 + \dots + N_L$ and $N_i = m_is_i$ so that $\mgg\gl_N$ has $\mgg\gl_{N_1} \oplus \dots \oplus \mgg\gl_{N_L}$ as subalgebras and we can take $f = f_1 + \dots + f_L$ with $f_i$ the rectangular element in the $i$-th summand. Assume that $0<s_1 < s_2 < \dots < s_L$, i.e. the partition is ordered. 
The first iteration is as Example \ref{exam2} with 
\[
M^{(1)} =N, \qquad M^{(1)} =  M^{(1)}_1 + M^{(1)}_2,\qquad 
 M^{(1)}_1 = N_1~, \qquad M^{(1)}_2 = N_2 + \dots + N_L.
 \]
  The reduction has an affine subVOA of type $\mgg\gl_{M^{(1)}_2}$ and the second iteration is the reduction on this subVOA again as in Example \ref{exam2} with 
  \[
  M^{(2)} = M^{(1)}_2, \qquad 
 M^{(2)} = M^{(2)} _1 + M^{(2)}_2, \qquad 
 M^{(2)} = M^{(1)}_2 ,\qquad  M^{(2)}_1 = N_2~, \qquad 
 M^{(2)}_2 = N_3 + \dots + N_L. 
  \]
  The reduction has an affine subVOA of type $\mgg\gl_{M^{(2)}_2}$ and one iterates with
  \[
   M^{(i)} = M^{(i-1)}_2 ,\qquad M^{(i)}= M^{(i)}_1 + M^{(i)}_2,\qquad  M^{(i)}_1 = N_i~, 
   \qquad 
  M^{(i)}_2 = N_{(i+1)} + \dots + N_L~. 
  \]
Then the complex for all iterations becomes
\[
V^k(\mgg\gl_N) \otimes F', \qquad F' = \bigotimes_{i=1}^L F^{m_i^2 s_i(s_i-1)} \otimes \bigoplus_{1 \leq  i < j \leq L} \widehat F^{2m_i m_j s_j(s_i-1)}.
\]
where 
\be\label{eq:ffalgC}
\widehat{F}^{2m_i m_j s_j(s_i-1)} := \begin{cases}  
F^{2 m_i m_j s_j(s_i-1)} & \quad r_i \ \text{odd} \\ 
F^{2 m_i m_j s_is_j} \otimes \beta\gamma^{m_i m_j s_j} & \quad s_i \ \text{even} \end{cases}
\ee
for $1\leq i< j \leq L$. 

As shown below, this yields the following identity 
\be\label{eq:mainformulaA} 
\begin{split}
 \text{ch}[\cW_L^k&(\mgg\gl_N, (\mgg\gl_{N_1}, \dots, \mgg\gl_{N_L}), (f_1, \dots, f_L))] 
 = \ch[\cW^k(\mgg\gl_N, f)] \prod_{1\leq i < j \leq L}  
\text{ch}[ \beta\gamma^{m_im_j (s_j-s_i)}],
\end{split}
\ee
where $\ch[\CA]$ denotes the character of the vertex algebra $\CA$. 
The step by step derivation of this identity proceeds as follows:
\begin{equation}\label{eq:mainformulaB} 
\begin{split}
 \text{ch}[\cW_L^k&(\mgg\gl_N, (\mgg\gl_{N_1}, \dots, \mgg\gl_{N_L}), (f_1, \dots, f_L))] \\
 &= \text{sch}[V^k(\mgg\gl_N) \otimes F'] \\
&= \ch[V^k(\mgg\gl_N)]  \prod_{i=1}^L \text{sch}[F^{m_i^2 s_i(s_i-1)}]   \prod_{1\leq i < j \leq L}  \text{sch}[\widehat F^{2m_i m_j s_j(s_i-1)}]  \\  
&= \ch[V^k(\mgg\gl_N)]  \prod_{i=1}^L
\text{sch}[F^{m_i^2s_i(s_i-1)}]   \prod_{1\leq i < j \leq L}  
\text{sch}[ F^{2m_i m_j s_j(s_i-1)}]  \\  
&= \ch[V^k(\mgg\gl_N)]  \prod_{i=1}^L
 \text{sch}[F^{m_i^2 s_i(s_i-1)}]   \prod_{1\leq i < j \leq L}  
 \text{sch}[ F^{2m_i m_j s_i(s_j-1)}] 
 \left(\text{sch}[ F^{2m_i m_j (s_j-s_i)}]\right)^{-1} \\  
&= \ch[V^k(\mgg\gl_N)]  \prod_{i=1}^L
\text{sch}[F^{m_i^2 s_i(s_i-1)}]   \prod_{1\leq i < j \leq L} 
 \text{sch}[ \widetilde F^{2m_i m_j s_i(s_j-1)}] 
\prod_{1\leq i < j \leq L}  \text{ch}[ \beta\gamma^{m_i m_j (s_j-s_i)}]\\  
&= \ch[\cW^k(\mgg\gl_N, f)] \prod_{1\leq i < j \leq L}  
\text{ch}[ \beta\gamma^{m_im_j (s_j-s_i)}]
\end{split}
\end{equation}
where $\widetilde F^{2m_i m_j s_i(s_j-1)}$ is defined in equation \eqref{eq:ffalgB} 
while $\widehat F^{2m_i m_j s_j(s_i-1)}$ is defined in equation \eqref{eq:ffalgC}. 

The first equality in \eqref{eq:mainformulaB} is that the character coincides with the supercharacter of the complex as all higher homologies vanish. 
The third equality follows as the supercharacters of
 $\widehat F^{2m_i m_j s_j(s_i-1)}$ and 
 $F^{2m_i m_j s_j(s_i-1)}$ coincide. 
The fifth equality uses that the supercharacters of 
$\widetilde F^{2m_i m_j s_i(s_j-1)}$ and $F^{2m_i m_j s_i(s_j-1)}$ coincide as well as that the $\beta\gamma$ character is the inverse of the supercharacter of two free fermions. Finally the last equality is the coincidence of the supercharacter of the complex with the $\cW$-algebra character. 

This coincidence of characters is expected to also hold on the level of actual VOAs, that is 
\begin{conj}\label{conjectureB}
$\cW_L^k(\mgg\gl_N, (\mgg\gl_{N_1}, \dots, \mgg\gl_{N_L}), (f_1, \dots, f_L))\cong \cW^k(\mgg\gl_N, f)\otimes  \bigotimes\limits_{1\leq i < j \leq L}  
\beta\gamma^{m_im_j(s_j-s_i)}$.
\end{conj} 

In order to conclude this section, note that a completely analogous identity holds for 
characters of universal Verma modules i.e. 
\be\label{eq:mainformulaC} 
\begin{split}
 \text{ch}^V[\cW_L^k&(\mgg\gl_N, (\mgg\gl_{N_1}, \dots, \mgg\gl_{N_L}), (f_1, \dots, f_L))] 
 = \ch^V[\cW^k(\mgg\gl_N, f)] \prod_{1\leq i < j \leq L}  
\text{ch}^V[ \beta\gamma^{m_im_j (s_j-s_i)}].
\end{split}
\ee

\subsection{Conformal weights and refined characters}\label{refcharsect} 

Let $V$ be a vertex algebra that is $\frac{1}{2}\mathbb Z$-graded by conformal weight. 
Let  $\{ F_pV \subset V | p \in \frac{1}{2}\mathbb Z\}$ be an increasing, exhaustive and separated filtration, that is 
\begin{enumerate}
\item $F_pV \subset F_{p + \frac{1}{2}}V$ 
\item $V = \bigcup\limits_p F_pV$
\item $0 = \bigcap\limits_p F_pV$
\end{enumerate}
Let $L_0$ be the Virasoro zero-mode of a conformal field of $V$. 
The filtration is compatible with the vertex algebra structure and the conformal structure if 
\begin{enumerate}
\item $| \, 0 \, \rangle \in F_0V\setminus F_{-\frac{1}{2}}V$
\item $T F_pV \subset F_pV$  ($T$ the derivation of $V$)
\item $L_0 F_pV \subset F_pV$ 
\item $v_n F_qV \subset F_{p+q}V$. for all $p, q \in \frac{1}{2} \mathbb Z$, $v \in F_pV$ and $n \in \frac{1}{2} \mathbb Z$. 
\end{enumerate}
The associated graded 
\[
\text{gr}^F V = \bigoplus_{p \in \frac{1}{2} \mathbb Z} F_pV/F_{p-\frac{1}{2}}V
\]
is then a vertex algebra, the associated graded of $V$ with respect to $F$. 
Assume that $V$ is graded by some Cartan subalgebra $\mfh$, such that
\[
h F_pV \subset F_pV 
\]
for $h$ in $\mfh$ then  $\mfh$ and $L_0$ induce an action on  $\text{gr}^F V$ and 
we can add an additional grading $g$ acting on $F_pV/F_{p-\frac{1}{2}}V$ by multiplication with $p$. 
 The refined character of $V$ is then 
\begin{equation}
 \text{ch}_{\rm ref}[V](y, z, q) = \text{tr}_{\text{gr}^F V } \left( q^{L_0 -\frac{c}{24}} z^{h} y^{-2g} \right). 
\end{equation}

\begin{exam}\label{betagammaex}
Let $V$ be the vertex algebra of $m_im_j(s_j-s_i)$ $\beta\gamma$-VOAs. It is strongly generated by 
$\{ \beta^t_{a,b} , \gamma^t_{a,b} | 1 \leq a \leq m_i, 1 \leq b \leq m_j, 1\leq t \leq s_j - s_i \}$.  
The OPE is 
\[
 \beta^t_{a,b} (z) \gamma^t_{a,b}(w)  = (z-w)^{-1}
\]
and all others are zero. 
Assigning the field $ \beta^t_{a,b} (z) $ (and hence any of its modes) the degree $t$ and the field  $ \gamma^t_{a,b} (z) $  the degree $1-t$ clearly gives a filtration such that the associated graded is commutative. 
Using the same variables as in equations \eqref{eq:charB} and \eqref{eq:charC}, the refined character is
\[
\bal
&{ \rm ch}_{\rm ref}[ \beta\gamma^{m_i m_j (s_j-s_i)}] = \\
 & \prod_{t=1}^{s_j-s_i} \ \ \ \ 
\prod_{m_1+\cdots + m_{i-1} +1 \leq a\leq 
m_1+\cdots+ m_i } \ \ \prod_{ m_1+\cdots + m_{j-1} +1 \leq b \leq 
m_1+\cdots+ m_j } (y^{-2t} z u_a u_b^{-1})^{-1}_\infty \
 (y^{2-2t} u_a^{-1}u_b)^{-1}_\infty~.\\
 \eal
 \]
The algebra of modes satisfies
\[
[\beta^t_{a,b, n}\, ,\, \gamma^t_{a,b, m}]  = \delta_{n+m, 0}
\]
and all other commutators vanish. Thus the $\beta^t_{a,b}$ form a commutative subalgebra and the submodule $B(m_i, m_j, s_j-s_i)$ of the $\beta\gamma$-VOA generated by this subalgebra acting on the vacuum has refined character \be\label{eq:Bcharacter}
\bal
 & {\rm ch}_{\rm ref}[ B(m_i, m_j, s_j-s_i)] = \\
&  \prod_{t=1}^{s_j-s_i} \ \ \ \ 
\prod_{m_1+\cdots + m_{i-1} +1 \leq a\leq 
m_1+\cdots+ m_i } \ \ \prod_{ m_1+\cdots + m_{j-1} +1 \leq b \leq 
m_1+\cdots+ m_j } (y^{-2t} z u_a u_b^{-1})^{-1}_\infty ~.
\eal
 \ee
\end{exam}

The standard Filtration for a vertex algebra is given by fixing a minimal set of strong generators and given $x_n$ degree $h$ if the strong generator $x(z)$ has conformal weight $h$ \cite{Li}. 
This depends on a choice of conformal vector. 
Recall that the strong generators of $\cW^k(\mgg, f)$ are labelled by $\alpha \in \Delta_{h.w.}$. If $\delta(\alpha) = j$, then the strong generator corresponds to a highest-weight vector of highest-weight $j$ with respect to the $\mathfrak{sl}_2$-triple for the quantum Hamiltonian reduction. It thus corresponds to the $u=2j+1$ dimensional representation. The conformal weigth of the generator is $j+1$ or in terms of the dimension $(u+1)/2$. 
From this we can read of that the refined character with respect to the standard filtration of Example \ref{exam} is 
\begin{equation}\label{eq:itcharA}
\text{ch}_{\rm ref}[V_\lambda]  = q^{-\frac{c}{24}} q^{\frac{(\lambda| \lambda + 2\rho) }{2(k+h^\vee)} -\lambda(x+y)}\ \prod_{i=1}^L X_i\ \prod_{1\leq i <j \leq L} X_{i,j}
\end{equation}
where 
\be\label{eq:itcharB}
\bal
X_i = & \prod_{t=1}^{s_i} (y^{-2t}z)_\infty^{-m_i}  \prod_{m_1+\cdots + m_{i-1} +1 \leq a<c\leq 
m_1+\cdots+ m_i } \ \ (y^{-2t} z u_au_c^{-1})^{-1}_\infty \ 
(y^{-2t} u_a^{-1} u_c)^{-1}_\infty\\
\eal 
\ee
and 
\be\label{eq:itcharC} 
\bal
X_{i,j} = & \prod_{t=1}^{s_i} \ \ \prod_{\substack{ m_1+\cdots + m_{i-1} +1 \leq a\leq 
m_1+\cdots+ m_i  \\ m_1+\cdots + m_{j-1} +1 \leq c \leq 
m_1+\cdots+ m_j }} \ \ (y^{-2t} y^{s_i- s_j}z u_au_c^{-1})^{-1}_\infty  
(y^{-2t} y^{s_i- s_j} u_a^{-1} u_c)^{-1}_\infty~.\\
\eal 
\ee
Let $w = \frac{1}{2} \sum_{k=1}^L s_kx_k$ and define a new conformal vector $L + \partial w$. The refined character for this filtration is 
\begin{equation}\label{eq:itcharD}
\text{ch}_{\rm ref}[V_\lambda]  = q^{-\frac{c}{24}} q^{\frac{(\lambda| \lambda + 2\rho) }{2(k+h^\vee)} -\lambda(x+y)}\ \prod_{i=1}^L X_i\ \prod_{1\leq i <j \leq L} X_{i,j} 
\end{equation}
where 
\be\label{eq:itcharE}
\bal
X_i = & \prod_{t=1}^{s_i}  (y^{-2t}z)_\infty^{-m_i} \prod_{m_1+\cdots + m_{i-1} +1 \leq a<c\leq 
m_1+\cdots+ m_i } \ \ (y^{-2t} z u_au_c^{-1})^{-1}_\infty \ 
(y^{-2t} u_a^{-1} u_c)^{-1}_\infty\\
\eal 
\ee
and 
\be\label{eq:itcharF} 
\bal
X_{i,j} = & \prod_{t=1}^{s_i} \ \ \prod_{\substack{ m_1+\cdots + m_{i-1} +1 \leq a\leq 
m_1+\cdots+ m_i  \\ m_1+\cdots + m_{j-1} +1 \leq c \leq 
m_1+\cdots+ m_j }} \ \ (y^{-2t} y^{2(s_i- s_j)}z u_au_c^{-1})^{-1}_\infty  (y^{-2t}  u_a^{-1} u_c)^{-1}_\infty~.\\
\eal 
\ee

In conclusion, note that the conformal weight refinement of identity \eqref{eq:mainformulaC} reads 
\be\label{eq:itcharG} 
\bal
 {\rm ch}^V_{\rm ref}[\cW_L^k&(\mgg\gl_N, (\mgg\gl_{N_1}, \dots, \mgg\gl_{N_L}), (f_1, \dots, f_L))] 
 = {\rm ch}_{\rm ref}^V[\cW^k(\mgg\gl_N, f)] \prod_{1\leq i < j \leq L}  
{\rm ch}_{\rm ref}[ \beta\gamma^{m_im_j (s_j-s_i)}].
\eal
\ee
where 
\be\label{eq:itcharH}
 {\rm ch}_{\rm ref}^V[\cW^k(\mgg\gl_N, f)]=\prod_{i=1}^L X_i \prod_{1\leq i <j \leq L} 
 X_{i,j}~. 
 \ee

\section{Comparison to partition functions of affine Laumon spaces}\label{comparison}

The goal of this section is to prove equation \eqref{eq:introG}, which relates the partition 
\eqref{eq:introC} to the refined universal Verma character 
associated to the $\beta$-truncation of the vertex algebra 
\[
\cW^k(\mgg\gl_N, f)\otimes  \bigotimes\limits_{1\leq i < j \leq L}  
\beta\gamma^{m_im_j(s_j-s_i)}.
\]
\begin{prop}\label{WZprop} 
Suppose the numerical invariants ${\bf r}=(r_0,\ldots, r_{\ell-1})$ are given by 
\be\label{eq:rformula} 
(r_0, \ldots, r_{\ell-1}) = 
\big( \underbrace{s_L, \ldots, s_L}_{m_L}\, , \, \underbrace{s_{L-1}, 
\ldots, s_{L-1}}_{m_{L-1}}, \ldots, \underbrace{s_1, \ldots, s_1}_{m_1}\big) 
\ee
for some $m_1,\ldots, m_L >0$ and $0< s_1 < \cdots < s_L$, $L\geq 1$. 
Then
\be\label{eq:PPtenF} 
\bal 
 Z_{\bf r} = {\rm ch}_{\rm ref}^V[\cW^k(\mgg\gl_N, f)]\prod_{1\leq i<j\leq L}\
{\rm ch}_{\rm ref}[B(m_i,m_j,s_j-s_i)].\\
\eal 
\ee
\end{prop}

{\it Proof}. 
As shown in Proposition \ref{newvarprop}, the partition function 
$Z_{\bf r}$ is given by 
\[
 \bal 
  Z_{\bf r}=\ & \prod_{a=0}^{\ell-1} \bigg(\prod_{t=1}^{r_a} \prod_{k=1}^\infty 
  {1\over 1-y^{-2t}z^k}  \bigg)  \\
  & 
\prod_{1\leq a< c\leq \ell}\ \bigg( 
\prod_{t=1}^{r_{\ell-c}}\prod_{k=1}^\infty {1\over 1- y^{-2t} z^{k} u_{a}u_{c}^{-1}}\bigg)\
\bigg( \prod_{t=1}^{r_{\ell-a}}\prod_{k=1}^\infty {1\over 1- y^{-2t} z^{k-1} u_{a}^{-1}u_{c}}\bigg). 
\eal 
 \]
The data $(m_1, \ldots, m_L)$ determines a partition 
\[
\{ 1, \ldots, \ell\} \simeq \bigcup_{i=1}^L \{m_1+\cdots + m_{i-1}+1, \ldots, 
m_1+ \cdots+ m_i\} 
\]
where the sets in the right hand side are pairwise disjoint. By convention, $m_{0}=0$. 
Moreover, given $1\leq c\leq \ell$,  one has 
\[
r_{\ell-c} = s_i 
\]
if and only if  $m_1+\cdots + m_{i-1}+1 \leq c \leq m_1+ \cdots+ m_i, \ 1\leq i \leq L$.  
Then the above expression for $Z_{\bf r}$ is equivalent to   
\be\label{eq:PPtenA}
Z_{\bf r} = \prod_{i=1}^L Z_i\ \prod_{1\leq i <j \leq L} Z_{i,j} 
\ee
where 
\be\label{eq:PPtenB} 
\bal
Z_i = & \prod_{t=1}^{s_i}  {(y^{-2t}z)_\infty^{m_i}}  \prod_{m_1+\cdots + m_{i-1} +1 \leq a<c\leq 
m_1+\cdots+ m_i } \ \ (y^{-2t} z u_au_c^{-1})^{-1}_\infty \ 
(y^{-2t} u_a^{-1} u_c)^{-1}_\infty\\
\eal 
\ee
and 
\be\label{eq:PPtenC} 
\bal
Z_{i,j} = & \prod_{t=1}^{s_j}  \prod_{v=1}^{s_i}
\ \ \prod_{m_1+\cdots + m_{i-1} +1 \leq a\leq 
m_1+\cdots+ m_i } \ \ \prod_{ m_1+\cdots + m_{j-1} +1 \leq c \leq 
m_1+\cdots+ m_j }
(y^{-2t} z u_au_c^{-1})^{-1}_\infty 
 (y^{-2v} u_a^{-1} u_c)^{-1}_\infty~.\\
\eal 
\ee
 Note that the $Z_i$ coincide with the factors $X_i$ defined in equation \eqref{eq:itcharE}. Furthermore, note that 
\be\label{eq:PPtenD}
\bal 
Z_{i,j}  & = \prod_{(a,c)} \prod_{t=1}^{s_j} \prod_{v=1}^{s_i} 
 \prod_{(a,c)}  (y^{-2t} z u_au_c^{-1})^{-1}_\infty
 (y^{-2v} u_a^{-1} u_c)^{-1}_\infty\\ 
 & =  \prod_{(a,c)} \prod_{t=1}^{s_j-s_i}(y^{-2t} z u_au_c^{-1})^{-1}_\infty\\ 
 & \ \ \ \prod_{(a,c)} \prod_{t=1}^{s_i} \prod_{v=1}^{s_i} 
  (y^{-2(t+s_j-s_i)} z u_au_c^{-1})^{-1}_\infty  (y^{-2v} u_a^{-1} u_c)^{-1}_\infty~.\\
\eal 
\ee
where $\prod_{(a,c)}$ stands for
\[
\prod_{m_1+\cdots + m_{i-1} +1 \leq a\leq 
m_1+\cdots+ m_i } \ \ \prod_{ m_1+\cdots + m_{j-1} +1 \leq c \leq 
m_1+\cdots+ m_j }~.
\]
By comparison to equations \eqref{eq:Bcharacter} and  \eqref{eq:itcharF}
one obtains 
\[
Z_{i,j} = X_{i,j}\, {\rm ch}_{\rm ref}[B(m_i,m_j,s_j-s_i)].
\]
In conclusion, using equation \eqref{eq:itcharH}, one obtains, 
\be\label{eq:PPtenE} 
\bal 
 Z  & =  \prod_{i=1}^L  X_i
\prod_{1\leq i<j\leq L}\ X_{i,j} \ {\rm ch}_{\rm ref}[B(m_i,m_j,s_j-s_i)] \\
& = {\rm ch}_{\rm ref}^V[\cW^k(\mgg\gl_N, f)]\prod_{1\leq i<j\leq L}\
{\rm ch}_{\rm ref}[B(m_i,m_j,s_j-s_i)].\\
\eal 
\ee

\hfill $\Box$

 \appendix
 
 \section{Some combinatorial identities}\label{combid} 
 This section proves equations \eqref{eq:morsecountA} and \eqref{eq:morsecountB}. 
 
Given a Young diagram $\mu$ and an integer $-\ell+1\leq c\leq \ell-1$, let 
\[
\bal
{\sf N}^{\geq 0}_1(\mu; c) & = \{(i,j)\in \mu\,|\, \mu'_i-j \equiv c\ {\rm mod}\ \ell\, \}\\
{\sf N}^{>0}_1(\mu; c) &  = 
\{(i,j)\in \mu\,|\, \mu'_i-j> 0,\ \mu'_i-j \equiv c\ {\rm mod}\ \ell\, \}\\
{\sf N}^{\geq 0}_2(\mu; c) & = 
\{(i,j)\in \mu\, |\, j-1\equiv c\ {\rm mod}\ \ell\, \}
\eal
\]
Moreover, for fixed $1\leq i \leq {\rm col}(\mu)$ let 
\[
\bal
{\sf N}^{\geq 0}_1(i; c) & = \{1\leq j \leq \mu_i'\, |\, \mu'_i-j \equiv c\ {\rm mod}\ \ell\, \}\\
{\sf N}^{>0}_1(i; c) & = 
\{1\leq j \leq \mu_i'-1 \, |\, \mu'_i-j \equiv c\ {\rm mod}\ \ell\, \}\\
{\sf N}^{\geq 0}_2(i; c) & = 
\{1\leq j \leq \mu_i'\, |\, j-1\equiv c\ {\rm mod}\ \ell\, \}\\
\eal
\]
where ${\rm col}(\mu)$ denotes the number of columns of $\mu$ and $\mu_i'$ denotes the height of the $i$-th column. 

As in equation \eqref{eq:colorfct} a partition $\mu$ will be said to be $a$-colored if and only if the box $(i,j)\in \mu$ is marked by $a-j+1$ mod $\ell$. 
\begin{lemm}\label{morselemmA}
Let $n_{a-c}$ denote the number
of boxes of $\mu$ marked by  $a-c$ mod $\ell$.  Then 
\be\label{eq:fixedcountABC}
n_{a-c}=|{\sf N}_1^{\geq 0}(\mu; c)| 
\ee
for any $-\ell+1 \leq c\leq \ell-1$. 
\end{lemm} 

{\it Proof}. 
By definition, $ n_{a-c}= |{\sf N}_2^{\geq 0}(\mu; c)|$. Hence 
it suffices to show that 
\be\label{eq:fixedcountAB}
|{\sf N}_1^{\geq 0}(\mu; c)| = |{\sf N}_2^{\geq 0}(\mu; c)|
\ee
for any $-\ell+1 \leq c\leq \ell-1$. 
Setting $j'=\mu_i'-j$, note that 
\[
\bal 
{\sf N}^{\geq 0}_1(i; c) & = \{1\leq j \leq \mu_i'\, |\, \mu'_i-j \equiv c\ {\rm mod}\ \ell\, \}\\
& = \{ 0\leq j' \leq \mu_i'-1 \, |\, j' \equiv c\ {\rm mod}\ \ell\, \},\\
\eal 
\]
which coincides with ${\sf N}^{\geq 0}_2(i; c)$ by definition. 
In conclusion, one obtains 
\[
|{\sf N}_1^{\geq 0}(i; c)| = |{\sf N}_2^{\geq 0}(i; c)|
\]
for any $1\leq i \leq {\rm col}(\mu)$, hence equation \eqref{eq:fixedcountAB} holds. 

\hfill $\Box$ 

\begin{lemm}\label{morselemmB} 
The following identity holds for any partition $\mu$: 
\be\label{eq:fixedcountD}
|{\sf N}_1^{> 0}(\mu; c)| = \left\{\begin{array}{ll}
|{\sf N}_2^{\geq 0}(\mu; c)| & {\rm for}\ -\ell+1\leq c\leq \ell-1, \ c\neq 0\\
& \\ 
|{\sf N}_2^{\geq 0}(\mu; c)|-{\rm col}(\mu) & {\rm for}\ c=0,\\
\end{array} \right.
\ee
where ${\rm col}(\mu)$ denotes the number of columns of $\mu$. 
\end{lemm}

{\it Proof}. 
The defining conditions of ${\sf N}_1^{> 0}(i; c)$ are equivalent to 
\[
1\leq j \leq \mu'_i-1,\qquad \mu'_i-j = k\ell + c, \ k \in \IZ.
\]
Suppose $\mu'_i = q_i \ell+ a_i$, $0\leq a_i\leq \ell-1$. Then these conditions are equivalent to 
\[
1\leq (q_i-k) \ell + a_i -c \leq q_i\ell+a_i-1 ~,
\]
hence also to 
\[
{1-c\over \ell} \leq k \leq q_i + {a_i-c-1\over \ell}~.
\]
Since $\ell \geq 2$, and $k,c,q_i,a_i$ are integers so that 
\[
-\ell+1\leq c\leq \ell-1, \qquad 0\leq a_i \leq \ell-1~, 
\]
these conditions are further  equivalent to 
\be\label{eq:fixedcountC}
\bal 
 0 \leq k \leq q_i -1\quad & {\rm if}\ 1\leq c\leq \ell-1\ {\rm and}\ 0\leq a_i\leq c \\
 0 \leq k \leq q_i \quad & {\rm if}\ 1\leq c\leq \ell-2\ {\rm and}\  c+1\leq a_i \leq \ell-1\\
  1 \leq k \leq q_i -1\quad & {\rm if}\ c=0\ {\rm and}\ a_i=0\\
  1 \leq k \leq q_i \quad & {\rm if}\ c= 0\ {\rm and}\ 1\leq a_i\leq \ell-1\\
 1 \leq k \leq q_i \quad & {\rm if}\ -\ell+1\leq c \leq -1\ {\rm and}\ 0\leq a_i \leq \ell+c\\
  1 \leq k \leq q_i+1 \quad & {\rm if}\ -\ell+1\leq c \leq -2\ {\rm and}\   \ell+c+1\leq a_i \leq \ell-1\\
 \eal
\ee
On the other hand, as shown in the proof of Lemma 
\ref{morselemmA}, the defining conditions of ${\sf N}_2^{\geq 0}(i; c)$ 
are 
\be\label{eq:fixedcountAB}
-{c\over \ell} \leq k \leq q_i + {a_i-c-1\over \ell}~.
\ee
These conditions are similarly equivalent to 
\be\label{eq:fixedcountB}
\bal 
 0 \leq k \leq q_i -1\quad & {\rm if}\ 0\leq c\leq \ell-1\ {\rm and}\ 0\leq a_i\leq c \\
 0 \leq k \leq q_i \quad & {\rm if}\ 0\leq c\leq \ell-2 \ {\rm and} \  c+1\leq a_i\leq \ell-1\\
 1\leq k \leq q_i \quad& {\rm if}\ -\ell+1\leq c\leq -1\ {\rm and}\ 0\leq a_i \leq \ell+c \\
 1\leq k \leq q_i+1 \quad& {\rm if}\ -\ell+1\leq c\leq -2\ {\rm and}\ \ell+c +1\leq a_i \leq \ell-1.\\
\eal
\ee
By comparison with \eqref{eq:fixedcountC}, one obtains the following table 
\[
\begin{array}{cccc} 
& & |{\sf N}_1^{> 0}(i; c)| & |{\sf N}_2^{\geq 0}(i; c)| \\ & & \\
1\leq c \leq \ell-1,\ 0\leq a_i\leq c   &  & q_i & q_i \\  & & \\
1\leq c\leq \ell-2,\ c+1\leq a_i \leq \ell-1 & & q_i+1 & q_i+1 \\  & & \\
c=0,\ a_i=0 & & q_i-1 & q_i \\   & & \\
c=0,\  1\leq a_i \leq \ell-1  & & q_i & q_i +1\\  & & \\
-\ell-1\leq c\leq -1,\ 0\leq a_i \leq \ell+c & & q_i & q_i \\  & & \\
-\ell-1\leq c\leq -2,\ \ell+c+1\leq a_i \leq \ell-1 & & q_i+1 & q_i+1 \\ & & \\
\end{array}
\]
In conclusion, 
\[
|{\sf N}_1^{> 0}(i; c)| = \left\{\begin{array}{ll}
|{\sf N}_2^{\geq 0}(i; c)| & {\rm for}\ -\ell+1\leq c\leq \ell-1, \ c\neq 0\\
& \\ 
|{\sf N}_2^{\geq 0}(i; c)|-1 & {\rm for}\ c=0.\\
\end{array} \right.
\]
This yields equation \eqref{eq:fixedcountD}. 

\hfill $\Box$

 \section{Some product identities}\label{infprodid} 
 
 This section proves equation \eqref{eq:PPfourDB}. 
 
  \begin{lemm}\label{partfctidA} 
 The following identity holds:
 \be\label{eq:PPfive}
\bal
& \prod_{\substack{1\leq a,c\leq \ell-1\\ a\neq c}} \prod_{t=1}^{r_a} 
{1\over 1- y^{-2t} z^{k}(\wq_{a-c})^{-1} (\wq_{a-c-1})^{-1} \cdots 
(\wq_{a-\ell+1})^{-1}} = \\
& \prod_{1\leq a<c\leq \ell-1} \bigg(\prod_{t=1}^{r_a} {1\over 1- y^{-2t} z^{k}(\wq_{a-c})^{-1} (\wq_{a-c-1})^{-1} \cdots 
(\wq_{a-\ell+1})^{-1}} \\
& \qquad \qquad \ 
\prod_{t=1}^{r_{a-c+\ell}}
{1\over 1- y^{-2t} z^{k-1}\wq_{a-c}\wq_{a-c-1} \cdots 
\wq_{a-\ell+1}}\bigg).
\eal
\ee
\end{lemm} 

{\it Proof}.
Note the set isomorphism 
\[
\{1\leq c < a\leq \ell-1\} \xrightarrow{\sim} \{ 1\leq a< c\leq \ell-1\} 
\]
\[
(a,c) \mapsto (a-c, \ell-c)
\]
Let $a'=a-c$ and $c' =\ell-c$. 
Since 
\[
z= \wq_a \wq_{a-1} \cdots \wq_{a-\ell+1}
\]
for any $a\in \IZ$, 
one has 
\[
\bal
& \prod_{1\leq c<a\leq \ell-1} \prod_{t=1}^{r_a} 
{1\over 1- y^{-2t} z^{k}(\wq_{a-c})^{-1} (\wq_{a-c-1})^{-1} \cdots 
(\wq_{a-\ell+1})^{-1}} = \\
& \prod_{1\leq c<a\leq \ell-1} \prod_{t=1}^{r_a} {1\over 1- y^{-2t} z^{k-1}\wq_{a}
\wq_{a-c-1} \cdots 
\wq_{a-c+1}} = \\
& \prod_{1\leq c<a\leq \ell-1} \prod_{t=1}^{r_a} {1\over 1- y^{-2t} z^{k-1}\wq_{a-\ell}
\wq_{a-c-1-\ell} \cdots 
\wq_{a-c+1-\ell}} = \\
&\prod_{1\leq a'<c'\leq \ell-1} \prod_{t=1}^{r_{a'-c'+\ell}} {1\over 1- y^{-2t} z^{k-1}
\wq_{a'-c'}\wq_{a'-c'-1} \cdots 
\wq_{a'-\ell+1}} = \\
&\prod_{1\leq a<c\leq \ell-1} \prod_{t=1}^{r_{a-c+\ell}} {1\over 1- y^{-2t} z^{k-1}\wq_{a-c}\wq_{a-c-1} \cdots 
\wq_{a-\ell+1}} \\
\eal
\]

\hfill $\Box$ 
 
 Recall that the variables \eqref{eq:newvar} are given by 
 \[
 u_a = (\wq_{-a} \cdots \wq \wq_{-\ell+1})^{-1},\ 1\leq a\leq \ell-1.
 \]
Then one has: 
\begin{lemm}\label{partfctidB}
Identity \eqref{eq:PPfive} can be further written as 
\be\label{eq:PPsix} 
\bal
& \prod_{\substack{1\leq a,c\leq \ell-1\\ a\neq c}} 
\prod_{t=1}^{r_a} {1\over 1- y^{-2t} z^{k}(\wq_{a-c})^{-1} (\wq_{a-c-1})^{-1} \cdots 
(\wq_{a-\ell+1})^{-1}} = \\
& \prod_{1\leq a< c\leq \ell-1}\
\bigg( 
\prod_{t=1}^{r_{\ell-c}}
{1\over 1- y^{-2t} z^{k} u_{a}u_{c}^{-1}}\bigg)\
\bigg( \prod_{t=1}^{r_{\ell-a}}{1\over 1- y^{-2t} z^{k-1} u_{a}^{-1}u_{c}}\bigg) \\
\eal
\ee
\end{lemm} 

{\it Proof}.
Again, note the set isomorphism 
\[
\{ 1\leq a< c\leq \ell-1\} \xrightarrow{\sim} \{ 1\leq a< c\leq \ell-1\}
\]
\[
(a,c) \mapsto (c-a, \ell-a) 
\]
Let $a' = c-a$ and  $c'= \ell-a$.
This yields 
\[
\bal
& \prod_{1\leq a<c\leq \ell-1} \bigg( 
\prod_{t=1}^{r_a} {1\over 1- y^{-2t} z^{k}(\wq_{a-c})^{-1} (\wq_{a-c-1})^{-1} \cdots 
(\wq_{a-\ell+1})^{-1}} \\
& \qquad \qquad \ \ \prod_{t=1}^{r_{a-c+\ell}} {1\over 1- y^{-2t} z^{k-1}\wq_{a-c}\wq_{a-c-1} \cdots 
\wq_{a-\ell+1}}\bigg) = \\
& \prod_{1\leq a'< c' \leq \ell-1}
\bigg( \prod_{t=1}^{r_{\ell-c'}} {1\over 1- y^{-2t} z^{k}(\wq_{-a'})^{-1} (\wq_{-a'-1})^{-1} \cdots 
(\wq_{-c'+1})^{-1}}\\
&\qquad \qquad \qquad \ \prod_{t=1}^{r_{\ell-a'}}
 {1\over 1- y^{-2t} z^{k-1}\wq_{-a'} \wq_{-a'-1} \cdots 
\wq_{-c'+1}}\bigg)=\\ 
& \prod_{1\leq a'< c' \leq \ell-1}\bigg( 
\prod_{t=1}^{r_{\ell-c'}}
{1\over 1- y^{-2t} z^{k} u_{a'}u_{c'}^{-1}}\bigg)\
\bigg( \prod_{t=1}^{r_{\ell-a'}}{1\over 1- y^{-2t} z^{k-1} u_{a'}^{-1}u_{c'}}\bigg) \\
\eal 
\]
as claimed in identity \eqref{eq:PPsix}. 

\hfill $\Box$

In conclusion, 
Lemmas \ref{partfctidA} and  \ref{partfctidB} yield identity \eqref{eq:PPfourDB} i.e.
\[
\bal
 & \prod_{a=1}^{\ell-1}\prod_{t=1}^{r_a}
\prod_{\substack{c=1\\ c\neq a}}^{\ell-1} 
{1\over (1-y^{-2t} z^k (\wq_{a-c})^{-1} (\wq_{a-c-1})^{-1} \cdots 
(\wq_{a-\ell+1})^{-1})} =\\
 & \prod_{1\leq a< c\leq \ell-1}\
\bigg( 
\prod_{t=1}^{r_{\ell-c}}
{1\over 1- y^{-2t} z^{k} u_{a}u_{c}^{-1}}\bigg)\
\bigg( \prod_{t=1}^{r_{\ell-a}}{1\over 1- y^{-2t} z^{k-1} u_{a}^{-1}u_{c}}\bigg). \\
\eal
\]

   \bibliography{AL_W_ref}
 \bibliographystyle{amsalpha}

\bigskip 
\addresses

\end{document}